\documentclass[11pt]{amsart}

\usepackage{amsmath,amssymb,amsthm,mathtools}
\usepackage{enumitem}
\usepackage{float}
\usepackage{microtype}
\usepackage{tikz}
\usetikzlibrary{arrows.meta,calc,shapes.geometric}
\definecolor{geomblue}{RGB}{35,83,120}
\tikzset{geompoint/.style={circle,fill=geomblue,inner sep=1.6pt},
 geomline/.style={draw=geomblue,line width=.7pt},
 geomlabel/.style={font=\small,inner sep=2pt},
 geomarrow/.style={geomline,-{Stealth[length=4pt]}}}
\usepackage{caption}
\usepackage[hidelinks]{hyperref}

\newcommand{\Hh}{\mathbb H}
\newcommand{\R}{\mathbb R}
\newcommand{\SL}{\mathrm{SL}}

\newcommand{\vp}{\boldsymbol v}
\newcommand{\grad}{\operatorname{grad}}

\theoremstyle{plain}
\newtheorem{theorem}{Theorem}[section]
\newtheorem{proposition}[theorem]{Proposition}
\newtheorem{lemma}[theorem]{Lemma}
\newtheorem{corollary}[theorem]{Corollary}
\theoremstyle{definition}
\newtheorem{definition}[theorem]{Definition}
\newtheorem{example}[theorem]{Example}
\newtheorem*{example*}{Example}
\theoremstyle{remark}
\newtheorem{remark}[theorem]{Remark}
\usepackage{etoolbox}
\AtEndEnvironment{remark}{\unskip\nobreak\hfill\(\diamond\)}

\title{Sharp Rigidity of Hyperbolic $p$-Barycenters}
\author[A. Elezi]{A. Elezi}
\address{Department of Mathematics and Statistics\\American University\\Washington, DC 20016}
\email{aelezi@american.edu}
\date{}

\begin{document}
\raggedbottom

\begin{abstract}
We study the hyperbolic $p$-barycenters $C_p$, $0\le p\le\infty$,
of $n$ points $\alpha_1,\ldots,\alpha_n$ in the hyperbolic plane.
For $0<p<\infty$, $C_p$ minimizes the sum of the $p$th powers of the
hyperbolic cosines of the distances to the points. At the endpoints,
$C_0$ minimizes the product of these cosines and $C_\infty$ their
maximum. The center map $p\mapsto C_p$ is real analytic at every
finite parameter. A point $P$ is the center for every parameter
exactly when, at each positive distance from $P$, the unit directions
to the points at that distance sum to zero. If the map is nonconstant,
the equation $C_p=P$ has at most $n-2$ solutions in $[0,\infty]$ for
each $P$, and at most $\lfloor(n-1)/2\rfloor$ when the points lie on a
geodesic. Finite solutions are counted with multiplicity, and the
endpoint contributes one if $C_\infty=P$. Each bound is attained in
arbitrarily small disks, with prescribed finite parameters and
multiplicities, and with or without the endpoint. The proofs
use exponential sums over distance shells, a variation-diminishing
property of the radial-weight kernel, and the Hobby--Rice theorem.
For the upper-half-plane roots of a totally complex real binary form,
$C_0$ is Julia's zero and $C_1$ is the hyperbolic zero. Their
coincidence forces the center map to be constant in degrees
at most six, but counterexamples exist in every even degree at least
eight. We obtain the same bounds on the number of distinct parameters
for Riemannian $L^q$ centers, the
minimizers of $\sum_i d(P,\alpha_i)^q$ on Hadamard manifolds
for $1<q<\infty$.
\end{abstract}

\subjclass[2020]{Primary 51M10; Secondary 11E76, 41A50, 53C20, 53C35.}
\keywords{Binary forms, Julia reduction, hyperbolic barycenter, circumcenter,
Riemannian center of mass, Riemannian $L^q$ centers, Hadamard manifolds, rigidity, exponential
sums, Descartes' rule of signs, variation diminishing kernels,
Hobby--Rice theorem.}

\maketitle

\section{Introduction}
Throughout, $\Hh^2$ is the hyperbolic plane of constant curvature $-1$
and $d$ is its distance. Configurations are nonempty finite multisets,
and all sums and cardinalities count multiplicity. For a configuration
$\mathcal A=\{\alpha_1,\ldots,\alpha_n\}$ in $\Hh^2$ and $0\le p<\infty$,
the hyperbolic $p$-center, or $p$-barycenter, $C_p(\mathcal A)$ is the
unique minimizer of the
\emph{$p$-potential}
\begin{equation}\label{eq:Phi-p}
 \Phi_{p,\mathcal A}(P)=\sum_{i=1}^n\varphi_p\bigl(d(P,\alpha_i)\bigr),
 \qquad
 \varphi_p(t)=
 \begin{cases}
  \dfrac{\cosh^pt-1}{p}, & p>0,\\[1ex]
  \log\cosh t, & p=0,
 \end{cases}
\end{equation}
and $C_\infty(\mathcal A)$ is the center of the smallest enclosing disk.
Since $\Phi_{0,\mathcal A}=\log\prod_i\cosh d(\cdot,\alpha_i)$ and
$\Phi_{1,\mathcal A}=\sum_i\cosh d(\cdot,\alpha_i)-n$, and since
$C_\infty(\mathcal A)$ minimizes $\max_i d(\cdot,\alpha_i)$, the points
$C_0$, $C_1$ and $C_\infty$ minimize, respectively, the geometric mean,
the arithmetic mean and the maximum of the numbers $\cosh d(P,\alpha_i)$. When the configuration is understood, we write
$C_p$ and $\Phi_p$.

This paper starts from the equilibrium laws and coincidence results of
\cite{EleziEquilibrium} for totally complex real binary forms. Up to a
constant factor, such a form is determined by the configuration
$\mathcal A$ of its upper-half-plane roots. Two reduction covariants,
Julia's zero and the hyperbolic zero, are characterized by
equilibrium laws of the form $\sum_i w(t_i)\,\vp_i=0$, with radial
weights $w=\tanh$ and $w=\sinh$, respectively
\cite[Propositions~9 and~4]{EleziEquilibrium}. Here $t_i$ is the
hyperbolic distance from the candidate center to the $i$-th
upper-half-plane root, and $\vp_i$ is the unit direction toward that
root; a root at the center contributes zero. For a sextic, the two zeros coincide precisely when its three
upper-half-plane roots, counted with multiplicity, form an equilateral
hyperbolic triangle or one of them is the hyperbolic midpoint of the
other two \cite[Proposition~12]{EleziEquilibrium}. In either case, the roots
at each positive distance from the common center balance on their own,
so changing the radial weight preserves equilibrium there.

The identity $\sinh r=\tanh r\,\cosh r$ places the two weights in the
family $w_p(r)=\tanh r\,\cosh^pr$. We show that, for finite $p$, $C_p$ is characterized by the equilibrium
law $\sum_i w_p(t_i)\,\vp_i=0$ (Proposition~\ref{prop:equilibrium});
hence $C_0$ is Julia's zero and $C_1$ is the hyperbolic zero (see also
Theorem~\ref{thm:endpoints}). The sextic phenomenon raises two questions for a general configuration.
When is the center map constant? If it is nonconstant, for how many
parameters can $C_p$ equal a given point $P$?

\subsection*{Main results}
If the center map of a configuration of $n$ points is
nonconstant, then for each $P$ the equation $C_p=P$ has at most $n-2$
solutions in $[0,\infty]$. When the points lie on a geodesic, the bound
improves to $\lfloor(n-1)/2\rfloor$. Finite solutions are counted with
multiplicity, and the endpoint contributes one if $C_\infty=P$
(Theorems~\ref{thm:n-minus-two} and~\ref{thm:collinear-multiplicity}).
Both bounds are sharp: configurations in arbitrarily small disks
attain them with the finite parameters and their multiplicities
prescribed, with or without the endpoint
(Theorem~\ref{thm:prescribed}).

Fix $P\in\Hh^2$ and group the points
of $\mathcal A$ other than $P$ by their distances $r_1<\cdots<r_s$ from
$P$ into \emph{shells}. Let the \emph{shell vector}
$\boldsymbol S_j(P)$ be the sum of the unit vectors at $P$ pointing to
the points of the $j$-th shell.

\emph{Shell criterion} (Theorem~\ref{thm:shell}). A point $P$ is a \emph{universal center}, meaning that $C_p=P$
for every $p\in[0,\infty]$, if and only if every shell vector
$\boldsymbol S_j(P)$ vanishes. It suffices that $C_p=P$ at $s$
distinct finite parameters. On a geodesic
this means reflection symmetry (Remark~\ref{rem:shell-geodesic}).

\emph{Center-multiplicity bound} (Theorems~\ref{thm:sign-bound}
and~\ref{thm:n-minus-two}). Suppose the center map is nonconstant. The
\emph{center multiplicity} $\kappa(P)$ is the sum of the
multiplicities of the finite solutions of $C_p=P$. For any unit vector $\boldsymbol e\in T_P\Hh^2$ for which the projections
$\langle\boldsymbol S_j(P),\boldsymbol e\rangle$ are not all zero,
$\kappa(P)$ is at most the number of sign changes in their sequence.
Moreover, if the $m$ nonzero shell vectors span
a $k$-dimensional subspace, then $\kappa(P)\le m-k$. A nonconstant three-point center map is injective on
$[0,\infty]$ and an immersion at finite parameters
(Corollary~\ref{cor:three-points}); both properties can fail for four
points.

\emph{Tail bound and points on a geodesic}
(Proposition~\ref{prop:hyp-tail} and Theorem~\ref{thm:collinear-multiplicity}).
For a unit vector $\boldsymbol e\in T_P\Hh^2$, let
\[
 N_{\boldsymbol e}(r)
 =\sum_{r_j>r}\langle\boldsymbol S_j(P),\boldsymbol e\rangle
\]
be the sum of the projected shell vectors beyond radius $r$. If $N_{\boldsymbol e}\not\equiv0$, the center multiplicity is at most
the number of sign changes of $N_{\boldsymbol e}$, which never exceeds the
sign bound. For points on a geodesic, all centers lie on it; when $P$
lies on that geodesic and $\boldsymbol e$ is tangent to it,
$N_{\boldsymbol e}$ is an integer walk, which gives the collinear bound
$\lfloor(n-1)/2\rfloor$ above.

Grouping the equilibrium law by shells expresses the net force
$F_p(P)=\sum_iw_p(t_i)\,\vp_i$ at $P$ as an exponential sum in $p$:
\[
 F_p(P)=\sum_{j=1}^s \tanh r_j\,e^{p\log\cosh r_j}\,\boldsymbol S_j(P).
\]
As $p$ increases, each shell force keeps its direction, while the
weights of more distant shells grow faster. Laguerre's rule of signs
therefore bounds the zeros of each nonzero scalar projection of $F_p(P)$
by the sign changes of the projected shell vectors.

A sharper bound comes from summing the shell vectors from the outside
in. Since $w_p(r)=\int_0^r\varphi_p''(t)\,dt$, the projected force
can also be written as
\[
 \langle F_p(P),\boldsymbol e\rangle
 =\int_0^\infty \varphi_p''(r)\,N_{\boldsymbol e}(r)\,dr.
\]
The variation-diminishing theorem for this kernel gives the tail bound. Combined with the Hobby--Rice theorem, it also gives the
configurations with prescribed parameters and multiplicities.

\subsection*{Binary forms}
For a totally complex real binary form $F$ with upper-half-plane
root configuration $\mathcal A$, write $\xi_p(F)=C_p(\mathcal A)$
for its \emph{$p$-zero}. Thus $\xi_0(F)$ is Julia's zero and
$\xi_1(F)$ is the hyperbolic zero. In degree $2n$, agreement at
$n-1$ distinct parameters in $[0,\infty]$ forces the $p$-zero
to be independent of $p$ (Corollary~\ref{cor:forms-n-minus-one}).
For quartics the family is always constant, and for sextics
agreement at any two distinct parameters suffices. In every even
degree at least $8$, however, any prescribed pair of distinct
parameters can give the same zero for a nonconstant family
(Subsection~\ref{subsec:thresholds}). Example~\ref{ex:integral-octic}
gives an integral octic for the pair $0,1$.
Theorem~\ref{prop:octic-classification} characterizes all octics
with a parameter-independent zero.

For collinear roots the threshold is $\lfloor(n-1)/2\rfloor+1$.
Hence for collinear octics any two distinct parameters again
suffice: the family is constant precisely when the four roots
can be paired with a common hyperbolic midpoint
(Corollary~\ref{prop:octic-universal}). The proofs here are
independent of \cite{EleziEquilibrium}.
Subsection~\ref{sec:reduction} records two consequences for reduction: a
universal center determines the reducing matrices independently of $p$,
and an interior reduction decision is stable under small changes of $p$,
with an explicit bound in terms of the velocity of the center map.

\subsection*{Riemannian \texorpdfstring{$L^q$}{Lq} centers}
On a Hadamard manifold $M$, the $L^q$ center $E_q$ minimizes
$\sum_i d(P,\alpha_i)^q$ for $1<q<\infty$; thus $E_2$ is the Riemannian
center of mass, and at $q=\infty$ we take $E_\infty$ to be the
circumcenter. The method carries over for two reasons: the radial
weight $r^{q-1}=e^{(q-1)\log r}$ is again exponential in the parameter,
which gives the shell criterion and the sign bound, and the integral
transform with kernel $r^{q-2}$ is variation diminishing, which gives
the tail bound. It follows that, for a nonconstant center path, the general bound
is $n-2$ and the collinear bound is $\lfloor(n-1)/2\rfloor$. Both are
attained at prescribed distinct finite parameters: the collinear bound
on any geodesic, and the general bound when $\dim M\ge2$
(Theorems~\ref{thm:Lq-rigidity} and~\ref{thm:Lq-collinear},
Corollary~\ref{thm:Lq-sharp}). Away from the data points, finite
solutions of $E_q=P$ can also be counted with multiplicity. At a data
point the force and the displacement of the center can vanish to
different orders, so the bounds there count distinct parameters
(Remark~\ref{rem:Lq-contact}).

\subsection*{Scope and related work}
Biswas~\cite{Biswas} introduced hyperbolic $p$-barycenters for
$1\le p\le\infty$ in proper geodesically complete CAT$(-1)$ spaces.
He used the associated barycenter maps to extend M\"obius
homeomorphisms of the boundaries and studied their limit as
$p\to\infty$. Here we fix a finite configuration in the hyperbolic
plane, include the range $0\le p<1$, and study how often the center
map takes a given value. The map is real analytic at every finite
parameter (Theorem~\ref{thm:path}). The counting arguments and
sharpness constructions also apply in $\Hh^d$ for $d\ge2$
(Remarks~\ref{rem:scope} and~\ref{rem:scope-constructions}).

For background on Riemannian centers of mass and $L^q$ centers, see
Karcher~\cite{Karcher} and Afsari~\cite{Afsari}. Da Silva Alves,
Oliveira and Greg\'orio~\cite{SilvaAlvesOliveiraGregorio} discuss
continuity of weighted $L^\alpha$ centers with respect to the exponent
for $1\le\alpha\le2$. On the real line, Blodgett~\cite{Blodgett} studied the path of
minimum $\ell^p$-norm location estimators as the exponent varies
from $1$ to $\infty$. Grouping observations at equal distances
from a candidate center, he applied Descartes' rule of signs
to the resulting exponential sum
\cite[Theorem~3.1 and Corollary~3.1]{Blodgett}.
He also used summation by parts to obtain a positivity criterion
from weighted tail sums \cite[Proposition~4.2]{Blodgett}.
Section~\ref{sec:Lq} applies the exponential balancing method
of Section~\ref{sec:rigidity} to $L^q$ centers on Hadamard
manifolds. On a geodesic, Theorem~\ref{thm:Lq-collinear}
gives the sharp bound $k=\lfloor(n-1)/2\rfloor$, including
the circumcenter endpoint, and attains it at any $k$ prescribed
distinct finite parameters. 

Sections~\ref{sec:prelim} and~\ref{sec:parameter} treat existence,
equilibrium and the regularity of the center map; Sections~\ref{sec:rigidity}
and~\ref{sec:sharp} prove the bounds and their sharpness;
Sections~\ref{sec:forms}, \ref{sec:Lq} and~\ref{sec:questions} treat
binary forms, Riemannian $L^q$ centers and further questions.

\section{Hyperbolic \texorpdfstring{$p$}{p}-barycenters and their
equilibrium laws}\label{sec:prelim}

\subsection{Existence and uniqueness for finite \texorpdfstring{$p$}{p}}

We use the upper-half-plane model $\Hh^2=\{x+iy:x\in\R,\ y>0\}$ with
metric $ds^2=(dx^2+dy^2)/y^2$. In this normalization,
\begin{equation}\label{eq:upper-distance}
 \cosh d(x+iy,u+iv)=\frac{(x-u)^2+y^2+v^2}{2yv}.
\end{equation}
Gradients, Hessians, inner products and norms of tangent vectors refer
to this metric. Let $\exp_P:T_P\Hh^2\to\Hh^2$ denote the exponential map
at $P$. For $P\ne A$, let $\vp_A(P)=\exp_P^{-1}(A)/d(P,A)$
denote the unit tangent vector at $P$ pointing toward $A$.
By the first variation of distance, the derivative of $d(\cdot,A)$
at $P$ in a unit direction $V$ is $-\langle V,\vp_A(P)\rangle$. Hence
\begin{equation}\label{eq:grad-distance}
 \grad d(\cdot,A)\big|_P=-\vp_A(P)\qquad(P\ne A).
\end{equation}
For every $p\ge0$, the radial derivatives are
\[
 \varphi_p'(r)=\tanh r\,\cosh^p r,\qquad
 \varphi_p''(r)=\cosh^{p-2}r\,(1+p\sinh^2r).
\]

\begin{lemma}\label{lem:strict-convex}\label{Hessian-positive}
For $A\in\Hh^2$, the function $f_A(P)=\log\cosh d(P,A)$ is real
analytic on $\Hh^2$. For every $p\ge0$,
\[
 \grad\bigl(\varphi_p(d(\cdot,A))\bigr)=e^{pf_A}\grad f_A.
\]
For $P\ne A$ we have
$\grad f_A(P)=-\tanh d(P,A)\,\vp_A(P)$, while $\grad f_A(A)=0$.
Moreover, $\varphi_p(d(\cdot,A))$ has positive-definite Hessian and
is strictly convex along every nonconstant geodesic.
\end{lemma}

\begin{proof}
By \eqref{eq:upper-distance}, $\cosh d(P,A)$ is a positive rational
function of the coordinates of $P$ with nonvanishing denominator.
Hence $f_A$ is real analytic, including at $A$, and
\[
 \varphi_p\bigl(d(\cdot,A)\bigr)=
 \begin{cases}
  \dfrac{e^{pf_A}-1}{p}, & p>0,\\[1ex]
  f_A, & p=0.
 \end{cases}
\]
The gradient identities now follow from the chain rule and
\eqref{eq:grad-distance}. At $A$ the gradient of $f_A$ vanishes
because $f_A$ has a minimum there.

For the Hessian, take a unit-speed geodesic $\gamma$ with $\gamma(0)=P$
and let $q(t)=\cosh d(\gamma(t),A)$. If $\rho=d(P,A)>0$, the hyperbolic
cosine law yields
\[
 q(t)=\cosh\rho\cosh t-\sinh\rho\cos\theta\sinh t,
\]
where $\theta$ is the angle between $\dot\gamma(0)$ and $\vp_A(P)$.
Hence $q''=q$ and $q^2-(q')^2=1+\sinh^2\rho\sin^2\theta\ge1$.
If $P=A$, then $q(t)=\cosh t$, so again $q''=q$ and
$q^2-(q')^2=1$. It follows that
\[
 \frac{d^2}{dt^2}\varphi_p\bigl(d(\gamma(t),A)\bigr)
 =q^{p-2}\bigl(q^2-(q')^2+p(q')^2\bigr)>0,
\]
also for $p=0$, when the left side is $(\log q)''$.
At $t=0$ this second derivative is the Hessian evaluated on
$\dot\gamma(0)$. Since the initial point and unit direction were
arbitrary, the Hessian is positive definite everywhere, and the
function is strictly convex along every nonconstant geodesic.
\end{proof}

\begin{theorem}\label{thm:exist-unique}
For every nonempty finite configuration $\mathcal A$ in $\Hh^2$ and
every $0\le p<\infty$, the potential $\Phi_{p,\mathcal A}$ is coercive
and strictly convex along geodesics. In particular, its minimizer
$C_p(\mathcal A)$ exists and is unique.
\end{theorem}

\begin{proof}
Strict convexity follows by summing Lemma~\ref{lem:strict-convex}.
As $P$ leaves every compact set, $d(P,\alpha_1)\to\infty$ and hence
$\varphi_p(d(P,\alpha_1))\to\infty$. The remaining summands are
nonnegative, so $\Phi_{p,\mathcal A}$ is coercive. Continuity and
compactness of closed bounded balls give existence; strict convexity
gives uniqueness.
\end{proof}

\subsection{Circumcenters and convex hulls}\label{subsec:circumcenters}

For $P\in\Hh^2$, the radius of
the smallest closed disk centered at $P$ that contains $\mathcal A$ is $r(P)=\max_i d(P,\alpha_i)$.
The function $r$ is continuous and coercive, so it attains its
minimum, the circumradius $R(\mathcal A)$. For uniqueness,
suppose that distinct points $P,Q$ both attain this minimum $R$.
Let $M$ be their midpoint and choose $i$ with $d(M,\alpha_i)=r(M)$.
Lemma~\ref{lem:strict-convex} with $p=1$ gives strict convexity of
$\cosh d(\cdot,\alpha_i)$, and therefore
\[
 \cosh r(M)
 <\tfrac12\bigl(\cosh d(P,\alpha_i)+\cosh d(Q,\alpha_i)\bigr)
 \le\cosh R.
\]
This contradicts the definition of $R$. The unique minimizer is the
circumcenter $C_\infty(\mathcal A)$.

\begin{proposition}[Convex-hull containment]\label{prop:convex-hull}
The closed hyperbolic convex hull $K$ of $\mathcal A$ contains
$C_p(\mathcal A)$ for every $0\le p\le\infty$. In particular, all the
centers lie in one compact convex set.
\end{proposition}

\begin{proof}
A closed hyperbolic ball containing $\mathcal A$ also contains $K$,
since the ball is convex. Hence $K$ is compact. Nearest-point
projection $\pi_K$ onto $K$ fixes each $\alpha_i$ and is
$1$-Lipschitz \cite[Proposition~II.2.4]{BridsonHaefliger}, so
$d(\pi_KP,\alpha_i)\le d(P,\alpha_i)$ for every $i$. It follows that projection
cannot increase the maximum distance $r$, and since $\varphi_p$ is
increasing it cannot increase $\Phi_{p,\mathcal A}$ either. By uniqueness of the
minimizer in each case, $\pi_K C_p=C_p$ for $0\le p\le \infty$.
\end{proof}

\begin{remark}[One and two points]\label{rem:two-points}
For a single point or two coincident points, every center is that point.
For two distinct points $\alpha_1,\alpha_2$, let $M$ be their midpoint
and write $d(\alpha_1,\alpha_2)=2h$. By
Proposition~\ref{prop:convex-hull}, every center lies on the segment
$[\alpha_1,\alpha_2]$. A point on this segment at distance $a$ from
$\alpha_1$ is at distance $2h-a$ from $\alpha_2$. For finite $p$,
strict convexity yields
\[
 \varphi_p(a)+\varphi_p(2h-a)\ge 2\varphi_p(h),
\]
with equality only at $a=h$. The maximum distance
$\max\{a,2h-a\}$ is likewise minimized at $a=h$.
It follows that $C_p=M$ for every $p\in[0,\infty]$, and a nonconstant center
map requires $n\ge3$.
\end{remark}

\subsection{Equilibrium and equivariance}

Fix $P\in\Hh^2$. As in the introduction, write $t_i=d(P,\alpha_i)$
and $\vp_i=\vp_{\alpha_i}(P)$ when $t_i>0$.

\begin{definition}[Forces]\label{def:forces}
The \emph{$p$-force exerted by $\alpha_i$ at $P$} is $w_p(t_i)\vp_i$,
where $w_p(t)=\varphi_p'(t)=\tanh t\,\cosh^p t$. Its value at $p=0$,
$\tanh(t_i)\vp_i$, is the \emph{base force}. When $P=\alpha_i$ the
force is zero, by continuous extension. The \emph{net $p$-force} is
\begin{equation}\label{eq:equilibrium-p}
 F_p(P)=\sum_{i=1}^n w_p(t_i)\vp_i .
\end{equation}
\end{definition}

\begin{proposition}[Equilibrium laws]\label{prop:equilibrium}
For $0\le p<\infty$, the net force is $F_p(P)=-\grad\Phi_p(P)$.
It follows that
\[
 P=C_p(\mathcal A)\quad\Longleftrightarrow\quad F_p(P)=0.
\]
For $p=\infty$, if all data points coincide then their common point is
the circumcenter. Otherwise, let $I(P)=\{i:t_i=r(P)\}$. Then
\begin{equation}\label{eq:circ-equilibrium}
 P=C_\infty(\mathcal A)\quad\Longleftrightarrow\quad
 0\in\operatorname{conv}\{\vp_i:i\in I(P)\},
\end{equation}
where the convex hull is taken in the tangent plane $T_P\Hh^2$.
In that case, for every unit vector $\boldsymbol u\in T_P\Hh^2$
there is an $i\in I(P)$ with
$\langle\boldsymbol u,\vp_i\rangle\le0$.
\end{proposition}

\begin{proof}
For finite $p$ and each $i$ with $t_i>0$, Lemma~\ref{lem:strict-convex}
yields
\[
 \grad\bigl(\varphi_p(d(\cdot,\alpha_i))\bigr)(P)
 =-\cosh^pt_i\,\tanh t_i\,\vp_i=-w_p(t_i)\,\vp_i,
\]
while this gradient is zero when $P=\alpha_i$. Summing yields
$\grad\Phi_p(P)=-F_p(P)$. By Theorem~\ref{thm:exist-unique}, $\Phi_p$ is strictly convex along geodesics and attains its minimum at $C_p$. A minimizer is a critical
point, and a strictly convex function has at most one critical point.
Hence $C_p$ is the unique zero of $\grad\Phi_p=-F_p$.

For the endpoint, if all data points coincide at $A$, then
$r(A)=0<d(Q,A)=r(Q)$ for every $Q\ne A$. Assume now that the data
contain two distinct points. Then $r(P)>0$, so all farthest
directions are defined. Suppose their convex hull misses zero, and
let $\boldsymbol u$ be its point closest to zero. For each
$i\in I(P)$ the segment from $\boldsymbol u$ to $\vp_i$ lies in the
hull, so the squared norm on it is minimized at $\boldsymbol u$;
hence $\langle\boldsymbol u,\vp_i-\boldsymbol u\rangle\ge0$ and
$\langle\boldsymbol u,\vp_i\rangle\ge\|\boldsymbol u\|^2>0$.
By \eqref{eq:grad-distance}, along $s\mapsto\exp_P(s\boldsymbol u)$
the derivative of the distance to $\alpha_i$ at $s=0$ is
$\langle\boldsymbol u,-\vp_i\rangle<0$ for all $i\in I(P)$. All
other distances are already strictly below $r(P)$. Since there are
only finitely many points, a sufficiently small move in this
direction decreases $r$, so $P$ is not the circumcenter.

Conversely, suppose zero is a convex combination of the farthest
directions. Given $Q\ne P$, write $Q=\exp_P(s\boldsymbol e)$ with
$s>0$ and $\boldsymbol e$ a unit vector. Taking the inner product
of the convex combination with $\boldsymbol e$ shows that some
$i\in I(P)$ satisfies $\langle\boldsymbol e,\vp_i\rangle\le0$;
since $\boldsymbol e$ is arbitrary, this proves the last assertion.
For this $i$, the hyperbolic cosine law yields
\[
 \cosh d(Q,\alpha_i)\ge\cosh s\,\cosh r(P)>\cosh r(P).
\]
Hence $r(Q)>r(P)$ for every $Q\ne P$, so $P$ is the circumcenter.
\end{proof}

\begin{proposition}[Isometry equivariance]\label{prop:equivariance}
For every isometry $g$ of $\Hh^2$ and every $0\le p\le\infty$,
$C_p(g\mathcal A)=gC_p(\mathcal A)$.
\end{proposition}

\begin{proof}
For finite $p$, distance invariance yields
$\Phi_{p,g\mathcal A}(gP)=\Phi_{p,\mathcal A}(P)$. The maximum distance
satisfies the same invariance. The result follows by uniqueness of
the minimizers.
\end{proof}

\subsection{Binary forms and \texorpdfstring{$p$}{p}-zero maps}\label{subsec:binary-forms}

Let $F$ be a totally complex real binary form of degree $2n\ge2$:
it is homogeneous and has no root in $\mathbb P^1(\R)$. In particular,
$F(1,0)\ne0$, and
\begin{equation}\label{eq:factor-form}
 F(X,Z)=a_0\prod_{i=1}^n(X-\alpha_iZ)(X-\overline{\alpha_i}Z),
 \qquad a_0\in\R\setminus\{0\},\quad \alpha_i\in\Hh^2.
\end{equation}
We define the \emph{$p$-zero} of $F$ by
$\xi_p(F)=C_p(\{\alpha_1,\ldots,\alpha_n\})$; for $p=\infty$ we also
call it the \emph{circumcenter zero}. This depends only on $F$, since
the multiset of upper-half-plane roots does not depend on the
labeling or on $a_0$. These
points are covariants by Proposition~\ref{prop:equivariance}. Indeed,
if $g=\left(\begin{smallmatrix}a&b\\c&d\end{smallmatrix}\right)
\in\SL_2(\R)$ and $F^g(X,Z)=F(aX+bZ,cX+dZ)$, the roots of $F^g$
are the images of the roots of $F$ under the inverse of
$g:z\mapsto(az+b)/(cz+d)$. Hence
$\xi_p(F^g)=g^{-1}\xi_p(F)$ for every $p$, including $p=\infty$.

Recall the two classical centers that motivated the family. A form of
degree $2n$ is stable if each root has multiplicity less than $n$.
For a totally complex real form this means precisely that the
upper-half-plane roots do not all coincide. For a stable form,
Julia's zero \cite{Julia,Cremona,StollCremona} is the unique
minimizer, over $t\in\R$ and $u>0$, of
$\widetilde F(t,u)/u^{2n}$ \cite[Proposition~5.1]{StollCremona},
where
\[
 \widetilde F(t,u)=|a_0|^2\prod_{i=1}^n
 (|t-\alpha_i|^2+u^2)(|t-\overline{\alpha_i}|^2+u^2).
\]
If all upper-half-plane roots coincide, $F$ is a scalar multiple of a
power of a positive definite real quadratic. In this case we take
Julia's zero to be the common upper-half-plane root, as in the
extension following Proposition~5.2 of \cite{StollCremona}.
The hyperbolic zero is the hyperbolic center of mass of the
upper-half-plane roots, equivalently the minimizer of
$\sum_i\cosh d(P,\alpha_i)$ \cite{Galperin,EleziShaska}.

\begin{theorem}\label{thm:endpoints}
For a totally complex real binary form $F$, with the conventions above,
$\xi_0(F)$ is Julia's zero and $\xi_1(F)$ is the hyperbolic zero.
\end{theorem}

\begin{proof}
Write $P=t+iu$ and $\alpha_i=x_i+iy_i$. The two factors of
$\widetilde F$ from a conjugate pair are equal, and by
\eqref{eq:upper-distance} each is $2uy_i\cosh d(P,\alpha_i)$. Hence
\[
 \frac{\widetilde F(t,u)}{u^{2n}}
 =|a_0|^2\!\left(\prod_{i=1}^n(2y_i)^2\right)
   \exp\!\bigl(2\Phi_{0,\mathcal A}(P)\bigr),
 \qquad \mathcal A=\{\alpha_1,\ldots,\alpha_n\}.
\]
The prefactor is positive and independent of $P$, so for a stable form
Julia's zero is $C_0(\mathcal A)$. If all upper-half-plane roots coincide at $A$, then Julia's zero is
$A$ by convention, and $C_0(\mathcal A)=A$ as well. Finally,
$\Phi_{1,\mathcal A}(P)=\sum_i\cosh d(P,\alpha_i)-n$, so
$C_1(\mathcal A)$ is the hyperbolic zero.
\end{proof}

\section{Regularity of the center map}\label{sec:parameter}

Fix the configuration $\mathcal A$. By Lemma~\ref{lem:strict-convex},
the Hessian of $\Phi_p$ is positive definite for every finite $p\ge0$.
We use this fact and the analytic implicit function theorem to show
that the center map $p\mapsto C_p$ is real analytic at every such parameter. We also
show that $C_p$ converges to the circumcenter $C_\infty$ as $p\to\infty$
by establishing an explicit upper bound for $d(C_p,C_\infty)$.
An analogous convergence for asymptotic hyperbolic $p$-barycenters has been shown by Biswas \cite[Definition~3.3, Proposition~3.7]{Biswas}. 

We use normal coordinates to describe the center map locally. Fix
$P\in\Hh^2$ and an orthonormal basis $(\boldsymbol e_1,\boldsymbol e_2)$
of $T_P\Hh^2$. The normal coordinates of $Q\in\Hh^2$ are
$x=(x_1,x_2)\in\R^2$, where
$\exp_P^{-1}(Q)=x_1\boldsymbol e_1+x_2\boldsymbol e_2$; we identify $x$
with this tangent vector. Since $\exp_P$ is a real-analytic
diffeomorphism, this is a global analytic chart centered at $P$. In the
proofs below, $P$ will be a center $C_{p_0}$. Let
\[f_i(x)=f_{\alpha_i}(\exp_P x)=\log\cosh d(\exp_P x,\alpha_i)\] with $f_{\alpha_i}$ as in
Lemma~\ref{lem:strict-convex}, and
\begin{equation}\label{eq:G-normal}
 \widetilde\Phi(p,x)=\sum_{i=1}^n\frac{e^{pf_i(x)}-1}{p},
 \qquad
 G(p,x)=\partial_x\widetilde\Phi(p,x)=\sum_{i=1}^ne^{pf_i(x)}\,\partial_xf_i(x),
\end{equation}
the first read as $\sum_if_i(x)$ at $p=0$, and $\partial_x$ denoting the
column of partial derivatives. Both functions are real analytic in
$(p,x)\in\R\times\R^2$. Notice that $\widetilde\Phi(p,x)=\Phi_p(\exp_Px)$ for
$p\ge0$. As in Section~\ref{sec:prelim}, let $t_i=d(P,\alpha_i)$ and
$\vp_i=\vp_{\alpha_i}(P)$.
\begin{enumerate}[label=\textup{(G\arabic*)},ref=\textup{(G\arabic*)},leftmargin=3.2em]
\item\label{G:zeros} For $p\ge0$, $G(p,x)=0$ if and only if
$\exp_Px=C_p$, since the minimizer is the only critical point of the
strictly convex function $\Phi_p$.
\item\label{G:origin} Since $f_i=f_{\alpha_i}\circ\exp_P$ and the differential
of $\exp_P$ at $0$ is the identity, the chain rule shows that
$\partial_xf_i(0)$ is the coordinate column of
$\grad f_{\alpha_i}(P)=-\tanh(t_i)\,\vp_i$
(Lemma~\ref{lem:strict-convex}). Since moreover
$e^{pf_i(0)}=\cosh^pt_i$, the definition \eqref{eq:equilibrium-p}
yields, as tangent vectors at $P$,
\begin{align*}
 G(p,0)&=-\sum_{i=1}^nw_p(t_i)\,\vp_i=-F_p(P),\\
 \partial_pG(p,0)&=-\sum_{i=1}^nw_p(t_i)\log\cosh(t_i)\,\vp_i ,
\end{align*}
a summand with $\alpha_i=P$ being zero. These identities hold for every
real $p$, where for $p<0$ we let $F_p(P)$ denote the sum
\eqref{eq:equilibrium-p}.
\item\label{G:hessian} The curves $t\mapsto\exp_P(tV)$ are geodesics
through $P$, so for $p\ge0$ and $V\in T_P\Hh^2$,
\[
 V^{\mathsf T}D_xG(p,0)\,V
 =\frac{d^2}{dt^2}\Phi_p\bigl(\exp_P(tV)\bigr)\Big|_{t=0}
 =\operatorname{Hess}_P\Phi_p(V,V).
\]
Hence $D_xG(p,0)$ is the matrix of the Hessian of $\Phi_p$ at $P$ in the
chosen basis, and it is positive definite by Lemma~\ref{lem:strict-convex}.
\end{enumerate}
The identification of $G$ with a force is made only at $x=0$; elsewhere
$G(p,x)$ is the column of partial derivatives of $\widetilde\Phi(p,\cdot)$.

\begin{theorem}[Regularity of the center map]\label{thm:path}
Fix a configuration $\mathcal A\subset\Hh^2$ and let $n=|\mathcal A|$.
\begin{enumerate}[label=\textup{(\roman*)}]
\item The center map $p\mapsto C_p(\mathcal A)$ is real analytic on $[0,\infty)$,
and at $p=0$ it extends analytically to a neighborhood of $0$. In
particular $C_p(\mathcal A)\to C_0(\mathcal A)$ as $p\downarrow0$.
\item For every $p>0$,
\begin{equation}\label{eq:quant-bound}
d(C_p,C_\infty)\leq\operatorname{arcosh}(n^{1/p}),
\end{equation}
hence $C_p(\mathcal A)\to C_\infty(\mathcal A)$ as $p\to\infty$.
\end{enumerate}
It follows that the center map
$\mathcal C_{\mathcal A}\colon[0,\infty]\to\Hh^2$,
$p\mapsto C_p(\mathcal A)$, is continuous.
\end{theorem}

\begin{proof}[Proof of \textup{(i)}]
Fix $p_0\ge0$ and let $P=C_{p_0}$. We use the normal coordinates centered at
$P$ and the map $G$ of \eqref{eq:G-normal}. Then $G(p_0,0)=0$ by
\ref{G:zeros}, and $D_xG(p_0,0)$ is positive definite, hence invertible,
by \ref{G:hessian}. The analytic implicit function theorem gives a real
analytic map $x(p)$ on an interval around $p_0$ with $x(p_0)=0$ and
\[
 G(p,x(p))=0 .
\]
For $p\ge0$ in this interval, \ref{G:zeros} gives $C_p=\exp_Px(p)$, so
the center map is real analytic near $p_0$. For $p_0=0$ the map $x(p)$ is
also defined for small $p<0$. There $\Phi_p$ and $C_p$ are not defined;
$\exp_Px(p)$ is a critical point of
$\sum_i(\cosh^pd(\cdot,\alpha_i)-1)/p$, and no minimizing property is
asserted. This extension is used only to make multiplicities at $p=0$
two-sided (Lemma~\ref{lem:order-contact}).
\end{proof}

\begin{proof}[Proof of \textup{(ii)}]
Let $R=r(C_\infty)=R(\mathcal A)$. If $R=0$, all points of
$\mathcal A$ coincide at $C_\infty$, so $C_p=C_\infty$ for every $p$ and there is
nothing to prove.

Suppose $R>0$. Minimality of $C_p$ yields
\begin{equation}\label{eq:radius-bound}
\cosh^p r(C_p)
\leq\sum_i\cosh^p d(C_p,\alpha_i)
\leq\sum_i\cosh^p d(C_\infty,\alpha_i)
\leq n\cosh^p R.
\end{equation}
For $Q\neq C_\infty$, let $\boldsymbol u$ be the unit tangent at $C_\infty$ toward $Q$.
By Proposition~\ref{prop:equilibrium}, some $\alpha_i\in \mathcal A$ of maximum distance from $C_\infty$ satisfies $\langle \boldsymbol u,\vp_{\alpha_i}(C_\infty)\rangle\le0$, that
is, the angle $\theta$ between $\boldsymbol u$ and $\vp_{\alpha_i}(C_\infty)$ is at least $\pi/2$.  The hyperbolic cosine law yields
\[
\cosh r(Q)\geq\cosh d(Q,\alpha_i)
\geq\cosh R\cosh d(C_\infty,Q).
\]
The same inequality holds when $Q=C_\infty$. Taking $Q=C_p$ and using
\eqref{eq:radius-bound} yields \[\cosh d(C_\infty,C_p)\leq n^{1/p},\]
which proves \eqref{eq:quant-bound}. 
\end{proof}

\begin{corollary}[Velocity of the center map]\label{cor:firstvariation}
Fix $0\le p<\infty$ and let $P=C_p$. Let $H_p$ be the Hessian of
$\Phi_p$ at $P$, viewed as an endomorphism of $T_P\Hh^2$ via the
hyperbolic metric. With $t_i=d(P,\alpha_i)$ and
$\vp_i=\vp_{\alpha_i}(P)$, the velocity of the center is given by
\begin{equation}\label{eq:velocity-center}
 \dot C_p
 = H_p^{-1}\sum_{i:\,t_i>0}
   w_p(t_i)\log\cosh(t_i)\,\vp_i .
\end{equation}
In particular, the center is stationary at $p$ if and only if the
sum on the right vanishes.
\end{corollary}

\begin{proof}
We prove the formula at an arbitrary $p_0\ge0$. Let $P=C_{p_0}$
and let $x(p)$ be the normal coordinates of $C_p$ centered at $P$, so
$x(p_0)=0$. By Theorem~\ref{thm:path}\textup{(i)}, $x$ is real analytic
near $p_0$, and by \ref{G:zeros}, $G(p,x(p))=0$. Differentiating at
$p_0$ yields
\[
 D_xG(p_0,0)\,\dot x(p_0)=-\partial_pG(p_0,0).
\]
By \ref{G:hessian}, the matrix $D_xG(p_0,0)$ represents $H_{p_0}$,
and since the differential of $\exp_P$ at $0$ is the identity,
$\dot x(p_0)$ represents $\dot C_{p_0}$ in the same basis.
On the right, the partial derivative keeps $P$ fixed, and
\ref{G:origin} yields
\[
 -\partial_pG(p_0,0)
 =\sum_{i:\,t_i>0}w_{p_0}(t_i)\log\cosh(t_i)\,\vp_i .
\]
Since $H_{p_0}$ is positive definite, it is invertible, and
\eqref{eq:velocity-center} follows.
\end{proof}

\begin{definition}[Multiplicity]\label{def:order-contact}
Let $p_0\ge0$ and $P=C_{p_0}$. Let $x(p)$ be the real analytic map
from the proof of Theorem~\ref{thm:path}\textup{(i)}, in normal
coordinates centered at $P$ with respect to an orthonormal basis of
$T_P\Hh^2$. Thus $x(p_0)=0$, $G(p,x(p))=0$, and $C_p=\exp_P x(p)$ for
$p\ge0$ near $p_0$. The map $x$ is also defined for small $p<0$ when
$p_0=0$. The \emph{multiplicity of $p_0$ as a solution of $C_p=P$},
denoted by $\operatorname{ord}_{p_0}(P)$, is the order of vanishing of
$x$ at $p_0$, that is, the least integer $k\ge1$ for which
\[
  \frac{d^k x}{dp^k}(p_0)\ne0 .
\]
If no such $k$ exists, we set $\operatorname{ord}_{p_0}(P)=\infty$. For
a point $Q\ne C_{p_0}$ we set $\operatorname{ord}_{p_0}(Q)=0$.
\end{definition}

The multiplicity does not depend on the basis. It is $1$ exactly when
$\dot C_{p_0}\ne0$; at higher multiplicity the center map is momentarily
stationary. Since $x$ is real analytic, $\operatorname{ord}_{p_0}(P)=\infty$
only if $C_p=P$ for every finite $p$.

\begin{lemma}\label{lem:order-contact}
Let $P\in\Hh^2$. At every finite parameter $p_0$ with $C_{p_0}=P$, the
map $p\mapsto F_p(P)$ of \eqref{eq:equilibrium-p} vanishes to order exactly
$\operatorname{ord}_{p_0}(P)$.
\end{lemma}

\begin{proof}
Let $G$ and $x(p)$ be as in Definition~\ref{def:order-contact}, with
normal coordinates centered at $P=C_{p_0}$. Then $G(p,x(p))=0$ for $p$
near $p_0$; for $p_0=0$ this includes small $p<0$. By \ref{G:origin}, $G(p,0)=-F_p(P)$,
also for $p<0$. Fix $p$ near $p_0$. By the chain rule,
\[
 \frac{d}{ds}G(p,sx(p))
 =D_xG(p,sx(p))\,x(p).
\]
Integrating over $0\le s\le1$  yields
\[
 F_p(P)=G(p,x(p))-G(p,0)=K(p)\,x(p), \quad
 K(p)=\int_0^1D_xG\bigl(p,s\,x(p)\bigr)\,ds .
\]
Since $x(p_0)=0$, $K(p_0)=D_xG(p_0,0)$ is the matrix of the Hessian of
$\Phi_{p_0}$ at $P$, by \ref{G:hessian}. The matrix $K(p)$ is real
analytic and $K(p_0)$ is positive definite. Hence $K(p)^{-1}$ is real
analytic near $p_0$, and multiplication by $K(p)$ preserves the order of
vanishing: $F_p(P)=K(p)x(p)$ and $x(p)=K(p)^{-1}F_p(P)$ vanish at $p_0$
to the same order, namely $\operatorname{ord}_{p_0}(P)$.
\end{proof}

\section{Universal centers and center multiplicity}\label{sec:rigidity}

In this section we ask when a point $P$ is the center of $\mathcal A$
for every parameter and, if it is not, for how many parameters
$C_p=P$. The first question has a geometric answer: $P$ is a universal
center if and only if, at each positive distance from $P$, the unit
directions from $P$ to the points of $\mathcal A$ at that distance sum
to zero. For a nonconstant center map, we then prove that the equation
$C_p=P$ has at most $n-2$ solutions in $[0,\infty]$, and at most
$\lfloor(n-1)/2\rfloor$ when the configuration lies on a geodesic.

\subsection{Universal centers}\label{subsec:universal}

Fix $P\in\Hh^2$, and let
\[
 0<r_1<\cdots<r_s
\]
be the distinct positive values among the distances $t_i=d(P,\alpha_i)$,
and let $\vp_i=\vp_{\alpha_i}(P)$ as in Section~\ref{sec:prelim}. They partition the points of $\mathcal A$ lying at positive
distances from $P$ into $s$ classes; we call the $j$-th class the $j$-th
\emph{shell} about $P$. The \emph{shell vector}, the
\emph{base shell force}, and the \emph{shell exponent} are defined respectively, by
\[
 \boldsymbol{S}_j(P):=\sum_{\{i:t_i=r_j\}}\vp_i,
 \qquad
 \boldsymbol b_j:=\tanh(r_j)\,\boldsymbol S_j(P),
 \qquad
 \lambda_j:=\log\cosh r_j .
\]
The vector $\boldsymbol b_j$ is the sum of the base forces exerted at $P$ by the
points of the $j$-th shell. Since $w_p(r)=\tanh r\,e^{p\log\cosh r}$, the
$p$-forces exerted at $P$ by these points add up to
$e^{p\lambda_j}\boldsymbol b_j$. Hence, grouped by shells, the net
$p$-force of Definition~\ref{def:forces} is
\begin{equation}\label{eq:shell-exp}
 F_p(P)=\sum_{j=1}^se^{p\lambda_j}\,\boldsymbol b_j .
\end{equation}
We call the $j$-th shell \emph{balanced} if $\boldsymbol S_j(P)=0$ and
\emph{unbalanced} otherwise. Balanced shells contribute nothing to
\eqref{eq:shell-exp} for any $p$. Varying $p$ rescales the shell forces without changing their directions.
Since $\lambda_1<\cdots<\lambda_s$, the weights of more distant shells
grow faster as $p$ increases.
\begin{lemma}[Zeros of exponential sums]\label{lem:exponential-zeros}
Let $\lambda_1<\cdots<\lambda_s$ be real numbers, and let
\[
 f(p)=\sum_{j=1}^s a_j e^{\lambda_j p}, \qquad a_j\in\R.
\]
If the coefficients are not all zero, the number of real zeros of $f$,
counted with multiplicity, is at most the number $V$ of sign changes in
$(a_1,\ldots,a_s)$ after zero entries are deleted. In particular, if $f$
vanishes at $s$ distinct real parameters, then $a_1=\cdots=a_s=0$.
\end{lemma}

This is Laguerre's extension of Descartes' rule of signs
\cite[Part~V, Chapter~1, Section~6, Problem~77]{PolyaSzego}. We include the proof because its Rolle argument, which counts zeros with multiplicity, is used again
in the proofs of Lemmas~\ref{lem:dirichlet} and~\ref{lem:laplace-zeros}.

\begin{proof}
Discard the zero coefficients; the exponents remain increasing. We
induct on $V$. If $V=0$, all terms of $f$ have the same sign and $f$ has
no zeros. If $V\ge1$, choose $k$ with $a_ka_{k+1}<0$ and
$\mu\in(\lambda_k,\lambda_{k+1})$. The function
$g(p)=e^{-\mu p}f(p)$ has the same zeros as $f$, and
\[
 g'(p)=\sum_j a_j(\lambda_j-\mu)e^{(\lambda_j-\mu)p}.
\]
The factor $\lambda_j-\mu$ is negative for $j\le k$ and positive for
$j>k$. Hence the coefficients of $g'$ are nonzero, they keep the sign
changes of $(a_j)$ within $\{1,\ldots,k\}$ and within $\{k+1,\ldots,s\}$,
and they lose the one between $k$ and $k+1$. By induction $g'$ has at
most $V-1$ zeros counted with multiplicity. If $g$ has distinct zeros
$z_1<\cdots<z_r$ of multiplicities $\mu_1,\ldots,\mu_r$, then $g'$ has a
zero of multiplicity $\mu_i-1$ at each $z_i$ and, by Rolle's theorem, a
further zero in each interval $(z_i,z_{i+1})$. Hence
$\sum_i\mu_i-1\le V-1$, and $g$, and so $f$, has at most $V$ zeros
counted with multiplicity.
\end{proof}

\begin{theorem}[Shell equilibrium criterion]\label{thm:shell}
Let $P\in\Hh^2$, with shells and shell vectors as above. The following are
equivalent:
\begin{enumerate}[label=\textup{(\roman*)}]
\item $P=C_p(\mathcal A)$ for every $p\in[0,\infty]$;
\item $P=C_{p_\ell}(\mathcal A)$ for $s$ distinct finite parameters
$p_1,\ldots,p_s\geq0$;
\item every shell about $P$ is balanced: $\boldsymbol{S}_j(P)=0$ for
$j=1,\ldots,s$.
\end{enumerate}
A point satisfying these conditions is called a \emph{universal center}.
When $s=0$, conditions (ii) and (iii) hold vacuously. 
\end{theorem}

\begin{proof}
If every point of $\mathcal A$ equals $P$, then $s=0$ and $C_p=P$ for
every $p\in[0,\infty]$, so (i) holds, and (ii) and (iii) hold
vacuously. Suppose $s\ge1$, that is, some point of $\mathcal A$ differs
from $P$.

(i)$\Rightarrow$(ii) is immediate.

(ii)$\Rightarrow$(iii): for every unit vector $\boldsymbol e\in
T_P\Hh^2$, the scalar exponential sum
$\langle F_p(P),\boldsymbol e\rangle=\sum_je^{\lambda_jp}
\langle\boldsymbol b_j,\boldsymbol e\rangle$ has $s$ distinct zeros, so
all its coefficients vanish by Lemma~\ref{lem:exponential-zeros}. Hence
every $\boldsymbol b_j$ is orthogonal to every $\boldsymbol e$, and
$\boldsymbol S_j(P)=0$ since $\tanh r_j>0$.

(iii)$\Rightarrow$(i): every shell force vanishes, so $F_p(P)=0$ and
$C_p=P$ for every finite $p$. The outermost
shell, at distance $r_s=\max_id(P,\alpha_i)$, is balanced,
$\sum_{i:\,t_i=r_s}\vp_i=0$, so zero lies in the
convex hull of the directions to the farthest points, and $P=C_\infty$ by
Proposition~\ref{prop:equilibrium}.
\end{proof}

\begin{remark}\label{rem:shell-geodesic}
Points equal to $P$ do not appear in the shell equilibrium condition because their
radial weight is zero for every $p$.  The theorem says that a universal center
is exactly a point at which every shell is directionally
balanced on its own.  In particular, any configuration invariant under a
nontrivial rotation about $P$ and decomposing into complete orbits on each
radius shell gives a universal center. For a configuration on a geodesic
$\ell$, shell balance at $P\in\ell$ means that, at each distance, the
points on the two sides of $P$ have equal multiplicities. Hence $P$ is a
universal center exactly when the configuration is invariant under the
reflection of $\ell$ fixing $P$. For instance, four points at signed
distances $\pm1$ and $\pm2$ from $P$ form two shells, each balanced on its
own.
\end{remark}

\begin{remark}[Power means and equidistant configurations]\label{prop:rho-monotone}
For $p>0$ let
\[
 M_{p,\mathcal A}(P)=\Bigl(\frac1n\sum_{i=1}^n\cosh^pd(P,\alpha_i)\Bigr)^{1/p}
 =\Bigl(1+\frac pn\Phi_{p,\mathcal A}(P)\Bigr)^{1/p}
\]
be the power mean of the numbers $\cosh d(P,\alpha_i)$, with limits
\[
 M_{0,\mathcal A}=\Bigl(\prod_i\cosh d(\cdot,\alpha_i)\Bigr)^{1/n},
 \qquad
 M_{\infty,\mathcal A}=\cosh\max_id(\cdot,\alpha_i).
\]
Then $M_{p,\mathcal A}$ is minimized at $C_p$ for every $p\in[0,\infty]$.
The minimum power mean is independent of $p$ exactly when all data
points coincide or all lie on a single balanced shell about some point
$P$; then $P$ is a universal center. Indeed, let $O=C_\infty$. If the
minimum is independent of $p$, the power-mean inequality at $O$ yields
\[
 \min M_{0,\mathcal A}\le M_{0,\mathcal A}(O)\le M_{\infty,\mathcal A}(O)
 =\min M_{\infty,\mathcal A}=\min M_{0,\mathcal A}.
\]
Hence $O=C_0$, and the geometric mean of the numbers
$\cosh d(O,\alpha_i)$ equals their maximum, so all distances from $O$
are equal. If they are positive, the points form a single shell about
$O$, which is balanced by Theorem~\ref{thm:shell}, since $C_0=O$.
Conversely, if the points lie on a single balanced shell of radius $r$
about $P$, then $P$ is a universal center and
$M_{p,\mathcal A}(P)=\cosh r$ for every $p$; if they coincide, the
minimum is $1$ for every $p$.
\end{remark}

\subsection{Center multiplicity}\label{subsec:counting}

By Theorem~\ref{thm:shell}, a point $P$ with $s$ shells is either a
universal center or satisfies $C_p=P$ for at most $s-1$ distinct
finite parameters. We now count each finite solution with its
multiplicity and use the directions of the shell vectors to sharpen the bound. 

Throughout this
subsection the center map is nonconstant, hence $n\ge3$
(Remark~\ref{rem:two-points}). Fix
$P\in\Hh^2$, with shells at the radii
$r_1<\cdots<r_s$, shell vectors $\boldsymbol S_j=\boldsymbol S_j(P)$, and
$\lambda_j=\log\cosh r_j$. For a unit vector $\boldsymbol e\in T_P\Hh^2$
let $V_{\boldsymbol e}(P)$ be the number of sign changes in the sequence
$\langle\boldsymbol S_1,\boldsymbol e\rangle,\ldots,
\langle\boldsymbol S_s,\boldsymbol e\rangle$, zero terms being ignored. Since balanced shells contribute nothing to $F_p(P)$, the unbalanced
shells determine the finite solutions of $C_p=P$
(Proposition~\ref{prop:equilibrium}) and, by
Lemma~\ref{lem:order-contact}, their multiplicities. Let $m\ge1$ be the number of unbalanced shells,
and $k\in\{1,2\}$ the dimension of the span of their shell vectors in $T_P\Hh^2$.

The center map is nonconstant, hence every multiplicity
$\operatorname{ord}_p(P)$ (Definition~\ref{def:order-contact}) is finite,
and $P$ is not a universal center, so by Theorem~\ref{thm:shell} it is
the center for at most $s-1$ finite parameters. The
\emph{center multiplicity} of $P$ is the finite sum
\[
 \kappa(P)=\sum_{0\le p<\infty,\ C_p=P}\operatorname{ord}_p(P).
\]
For a nonconstant center map, the \emph{fiber multiset over $P$}
consists of the finite parameters $p$ with $C_p=P$, each counted
with its multiplicity $\operatorname{ord}_p(P)$, together with
one occurrence of $\infty$ if $C_\infty=P$. Its finite part has
cardinality $\kappa(P)$, and its total cardinality is
$\kappa(P)+\mathbf 1_{\{C_\infty=P\}}$. 

\begin{theorem}[Sign and shell-span bounds]\label{thm:sign-bound}
\begin{enumerate}[label=\textup{(\roman*)}]
\item If the numbers $\langle\boldsymbol S_j,\boldsymbol e\rangle$ are not
all zero, then $\kappa(P)\le V_{\boldsymbol e}(P)$, the number of sign
changes in this sequence. In particular, if they all
have the same sign, $P$ is not the $p$-center for any finite $p$.
\item \textup{(Shell-span bound.)} If the shell vectors of the $m$
unbalanced shells span a $k$-dimensional subspace, then
$\kappa(P)\le m-k$.
\item If $m\ge2$ and $\kappa(P)=m-1$, then the unbalanced shell vectors
are parallel and alternate in direction when ordered by increasing
radius. More explicitly, if $j_1<\cdots<j_m$ index the unbalanced shells,
then
\[
 \boldsymbol S_{j_\ell}=(-1)^{\ell-1}c_\ell\,\boldsymbol u,\qquad c_\ell>0,
\]
for a unit vector $\boldsymbol u$.
\end{enumerate}
\end{theorem}

\begin{proof}
(i) The scalar function $\langle F_p(P),\boldsymbol e\rangle$ is a
nonzero exponential sum whose coefficients
$\tanh(r_j)\langle\boldsymbol S_j,\boldsymbol e\rangle$ have the signs of
the numbers $\langle\boldsymbol S_j,\boldsymbol e\rangle$. By
Lemma~\ref{lem:order-contact} it vanishes at each finite parameter of the
fiber to order at least its multiplicity, and Lemma~\ref{lem:exponential-zeros} gives
$\kappa(P)\le V_{\boldsymbol e}(P)$.

(ii) Since balanced shells contribute zero terms, we may omit them and
index the unbalanced shells by $1,\ldots,m$ in order of distance. If
$k=1$, their shell vectors are parallel to a unit vector
$\boldsymbol u$. The sequence $\langle\boldsymbol S_j,\boldsymbol u\rangle$
has exactly $m$ nonzero terms, so $V_{\boldsymbol u}(P)\le m-1$, and (i)
gives $\kappa(P)\le m-1$. If $k=2$, choose a unit vector $\boldsymbol e\perp\boldsymbol
S_m$. Since the shell vectors span $T_P\Hh^2$, the numbers
$\langle\boldsymbol S_j,\boldsymbol e\rangle$, $1\le j\le m-1$, are not
all zero. As $\langle\boldsymbol S_m,\boldsymbol e\rangle=0$, at most
$m-1$ of the $m$ numbers are nonzero, so $V_{\boldsymbol e}(P)\le m-2$,
and (i) gives
$\kappa(P)\le m-2$. Thus a second independent direction allows a
projection that removes one shell, and lowers the bound by one.

(iii) By (ii), $m-1=\kappa(P)\le m-k$, so $k=1$. Choose a unit vector
$\boldsymbol u$ parallel to the unbalanced shell vectors, oriented so
that the first has positive projection. By (i), their $m$ nonzero
projections $\langle\boldsymbol S_{j_\ell},\boldsymbol u\rangle$ have at
least $m-1$ sign changes, so they alternate in sign; the balanced shells
contribute zero terms and are omitted.
\end{proof}

If $m=n$, then no point equals $P$ and every shell consists of a single
point. Hence, if $\kappa(P)=n-1$, part (iii) forces the points to lie on a
geodesic through $P$, on alternating sides in order of distance. The next
proof rules out equilibrium in this exceptional configuration.

\begin{theorem}[Center-multiplicity bound]\label{thm:n-minus-two}
Let $\mathcal A$ consist of $n\ge3$ points, counted with multiplicity, and
suppose that the center map $\mathcal C_{\mathcal A}$ is nonconstant. Then for
every $P\in\Hh^2$,
\[
 \kappa(P)+\mathbf 1_{\{C_\infty=P\}}\le n-2 .
\]
Here $\mathbf 1_{\{C_\infty=P\}}$ equals $1$ if $P=C_\infty$ and $0$
otherwise; the endpoint parameter is counted once. In particular, every fiber of a nonconstant center map has total
multiplicity at most $n-2$ on $[0,\infty]$, and if $C_p(\mathcal A)$ takes the same value at
$n-1$ distinct parameters in $[0,\infty]$, the center map is constant.
\end{theorem}

\begin{proof}
With $m$ and $k$ as above, $\kappa(P)\le m-k$ by
Theorem~\ref{thm:sign-bound}(ii).

Suppose first that $C_\infty\neq P$. Since $m\le n$ and $k\ge1$, we have
$\kappa(P)\le n-2$ unless $k=1$, $m=n$ and $\kappa(P)=n-1$. In that case no point equals $P$,
every shell is a single point, $s=n$, and the points lie on one geodesic
through $P$. By Theorem~\ref{thm:sign-bound}(iii), after labeling the points so
that $0<r_1<\cdots<r_n$, we have $\vp_j=(-1)^{j-1}\boldsymbol u$, i.e. the
points alternate sides of $P$ in order of distance. At a parameter $p_1$
with $C_{p_1}=P$ the
equilibrium law reads $\sum_{j=1}^n(-1)^{j-1}c_j=0$ with
$c_j=\tanh(r_j)\cosh^{p_1}(r_j)$. Since $p_1\ge0$, the radial weight
increases with distance, $0<c_1<\cdots<c_n$, and grouping the terms in pairs from the outside in we obtain
\[
 \sum_{j=1}^n(-1)^{j-1}c_j=(-1)^{n-1}\left[(c_n-c_{n-1})+(c_{n-2}-c_{n-3})+\cdots\right],
\]
with a final term $c_1$ when $n$ is odd. The bracket is positive, so the
sum is nonzero.
This contradiction yields $\kappa(P)\le n-2$.

Suppose now that $C_\infty=P$, and let the outermost shell, at distance
$r_s=\max_id(P,\alpha_i)$, contain $b$ points, counted with multiplicity.
By Proposition~\ref{prop:equilibrium}, zero lies in the convex hull of
their directions, so $b\ge2$. If $b=2$, the two directions are opposite,
so $\boldsymbol S_s=0$, and since $s\le n-b+1=n-1$ we get $m\le s-1\le
n-2$. If $b\ge3$, then $m\le s\le n-b+1\le n-2$. In both cases
$\kappa(P)\le m-k\le n-3$, and adding the endpoint $p=\infty$ gives the bound
$n-2$.
\end{proof}

The bound is sharp in the following sense: every multiset of
cardinality $n-2$ in $[0,\infty]$, with $\infty$ occurring at most
once, is the fiber multiset over $P$ of a configuration of $n$
points in any prescribed neighborhood of $P$
(Theorem~\ref{thm:prescribed}\textup{(ii)}).

\begin{corollary}[Three points]\label{cor:three-points}
For three points, counted with multiplicity, the center map is either constant,
or injective on $[0,\infty]$ and an immersion on $[0,\infty)$.
\end{corollary}

\begin{proof}
Suppose the center map is nonconstant. For $n=3$,
Theorem~\ref{thm:n-minus-two} gives at most one solution of $C_p=P$
in $[0,\infty]$ for each $P$, so the map is injective. Every finite
solution has multiplicity one, hence $\dot C_p\ne0$ for every finite $p$.
\end{proof}

\noindent Both conclusions fail for every $n\ge4$: by
Theorem~\ref{thm:prescribed}(ii), the center map can take the same value
at two prescribed parameters, or have zero velocity at a prescribed
finite parameter.

\subsection{Tail counts and collinear configurations}\label{subsec:tail}

\noindent The sign bound treats the shell vectors separately. Since
$w_p(r)$ increases with $r$ for $p\ge0$, we can instead assemble the net
force by successive positive increments: the first acts on all shells,
and each subsequent increment acts only on the shells farther out. This
suggests considering the sums of shell vectors beyond each radius. 

For $r\ge0$ define the \emph{tail vector}
\[
 \boldsymbol T_P(r):=\sum_{j:\,r_j>r}\boldsymbol S_j(P)
 =\sum_{i:\,t_i>r}\vp_i\in T_P\Hh^2,
\]
the sum of the shell vectors of the shells farther than $r$ from $P$,
equivalently of the unit directions toward the points farther than $r$.
Balanced shells contribute nothing to it. Let $r_0=0$. The tail vector
is constant on each interval $[r_{j-1},r_j)$ and vanishes on
$[r_s,\infty)$. As $r$ decreases past a shell radius $r_j$ the tail vector jumps
by the shell vector $\boldsymbol S_j(P)$. Its successive values
\[
 0,\quad \boldsymbol S_s,\quad \boldsymbol S_s+\boldsymbol S_{s-1},\quad
 \ldots,\quad \boldsymbol S_s+\cdots+\boldsymbol S_1
\]
are the vertices of a polygonal path in $T_P\Hh^2$, the
\emph{tail path}: it starts at the origin, is traversed from the
outer shells inward, and its edges are the shell vectors.

Using $\boldsymbol S_j(P)=\boldsymbol T_P(r_{j-1})-\boldsymbol T_P(r_j)$,
$w_p(0)=0$ and $\boldsymbol T_P(r_s)=0$, summation by parts yields
\begin{equation}\label{eq:abel}
 F_p(P)=\sum_{j=1}^s w_p(r_j)\,\boldsymbol S_j(P)
 =\sum_{j=0}^{s-1}\bigl(w_p(r_{j+1})-w_p(r_j)\bigr)\,
 \boldsymbol T_P(r_j).
\end{equation}
The weight $w_p$ increases with $r$, so $F_p(P)$
is a positive combination of the vertices of the tail path. Since
$w_p(r_{j+1})-w_p(r_j)=\int_{r_j}^{r_{j+1}}\varphi_p''(r)\,dr$, we
obtain the \emph{tail identity} in $T_P\Hh^2$:
\begin{equation}\label{eq:hyp-tail}
 F_p(P)=\int_0^\infty\varphi_p''(r)\,\boldsymbol T_P(r)\,dr.
\end{equation}
In this representation, the dependence on $p$ lies
entirely in the kernel $\varphi_p''$, while the configuration and the
base point $P$ determine the tail vector function. For a unit vector
$\boldsymbol e\in T_P\Hh^2$, the \emph{tail count}
\[
 N_{\boldsymbol e}(r)=\langle\boldsymbol T_P(r),\boldsymbol e\rangle
 =\sum_{j:\,r_j>r}\langle\boldsymbol S_j(P),\boldsymbol e\rangle
\]
is the projection of the tail vector, and
$\langle F_p(P),\boldsymbol e\rangle=\int_0^\infty\varphi_p''(r)\,
N_{\boldsymbol e}(r)\,dr$. The tail count $N_{\boldsymbol e}$ is real-valued. Like $\boldsymbol T_P$, it is constant on each interval $[r_{j-1},r_j)$ and vanishes on
$[r_s,\infty)$.

More generally, a \emph{finite step function} on $(0,\infty)$ is a
real-valued function that is constant on each of finitely many open
intervals covering $(0,\infty)$ up to finitely many breakpoints. For such $h$, let $S^-(h)$ be the number of sign changes in the sequence of its values on these open
intervals, taken in order, omitting intervals on which $h=0$ and ignoring
the values at the breakpoints. This is Karlin's notation $S^-$ from total positivity \cite{Karlin}, applied to the sequence of values of $h$ on its intervals. 

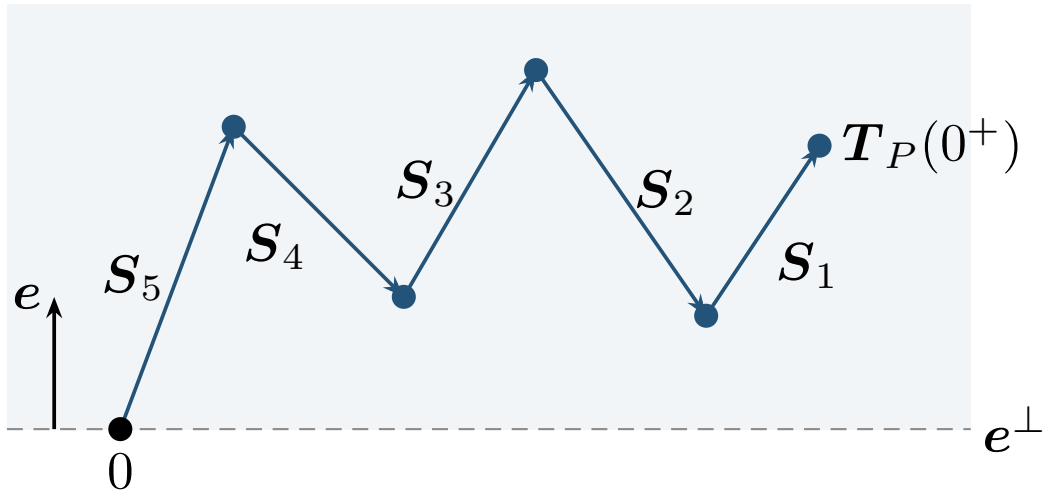
\begin{figure}[tbp]
\centering
\begin{tikzpicture}[every node/.style={geomlabel},x=1.25cm,y=1.25cm]
 \fill[geomblue!6] (-0.6,0) rectangle (4.5,2.25);
 \draw[black!45,densely dashed] (-0.6,0)--(4.5,0)
   node[right,black] {$\boldsymbol e^\perp$};
 \draw[geomarrow,black] (-0.35,0)--(-0.35,0.7) node[left] {$\boldsymbol e$};
 \coordinate (T0) at (0,0);
 \coordinate (T1) at (0.6,1.6);
 \coordinate (T2) at (1.5,0.7);
 \coordinate (T3) at (2.2,1.9);
 \coordinate (T4) at (3.1,0.6);
 \coordinate (T5) at (3.7,1.5);
 \draw[geomarrow] (T0)--(T1) node[midway,left] {$\boldsymbol S_5$};
 \draw[geomarrow] (T1)--(T2) node[midway,below left] {$\boldsymbol S_4$};
 \draw[geomarrow] (T2)--(T3) node[midway,left] {$\boldsymbol S_3$};
 \draw[geomarrow] (T3)--(T4) node[midway,right] {$\boldsymbol S_2$};
 \draw[geomarrow] (T4)--(T5) node[midway,below right] {$\boldsymbol S_1$};
 \foreach \p in {T1,T2,T3,T4,T5} {\node[geompoint] at (\p) {};}
 \node[geompoint,fill=black,label=below:{$0$}] at (T0) {};
 \node[right=2pt] at (T5) {$\boldsymbol T_P(0^+)$};
\end{tikzpicture}
\caption{A tail path with $s=5$ shells. Its vertices are the
values of $\boldsymbol T_P(r)$, starting from $0$ for $r$ beyond the
farthest point and adding the shell vectors $\boldsymbol S_5,\ldots,
\boldsymbol S_1$ from the outside in. The projections of the edges onto
$\boldsymbol e$ alternate in sign four times, so
$V_{\boldsymbol e}(P)=4$, while all vertices after the first lie in the
same open half-plane, so $S^-(N_{\boldsymbol e})=0$ (schematic)}
\label{fig:tail}
\end{figure}

The values of $N_{\boldsymbol e}$ are the projections onto
$\boldsymbol e$ of the vertices of the tail path, so
$S^-(N_{\boldsymbol e})$ counts how often these vertices alternate
between the two open half-planes bounded by $\boldsymbol e^\perp$,
vertices on $\boldsymbol e^\perp$ being omitted. On the other hand,
$V_{\boldsymbol e}(P)$ counts how often the edges of the tail path
alternate in direction across $\boldsymbol e^\perp$. Because the path
starts at the origin, every alternation between nonzero vertex
projections requires a corresponding alternation in the signs of the
projected edges. Hence vertex-side alternations cannot outnumber
edge-direction alternations,
$S^-(N_{\boldsymbol e})\le V_{\boldsymbol e}(P)$
(Proposition~\ref{prop:hyp-tail}). The inequality can be strict, since
the edges may zigzag while the vertices stay on one side
(Figure~\ref{fig:tail}).

The main result of this subsection is that the kernel $\varphi_p''$ is
variation diminishing for $p\ge0$: for a finite step function $h$ of
bounded support that is not zero almost everywhere, the transform
\[
 p\longmapsto\int_0^\infty\varphi_p''(r)h(r)\,dr
\]
has no more zeros, counted with multiplicity, than $h$ has sign changes.
This is a continuous analogue of Lemma~\ref{lem:exponential-zeros}.
There the discrete kernel $e^{p\lambda_j}$ acts on the coefficient
sequence $(a_j)$ by summation; here the radial kernel $\varphi_p''(r)$
acts on the step function $h$ by integration. In both cases, the
number of zeros, counted with multiplicity, is bounded by the number
of sign changes of the input. For the connection between variation
diminution and total positivity, see Karlin~\cite{Karlin}. The proof of the general result is given in
Subsection~\ref{subsec:proof-vd} and can be skipped on a first reading.

\begin{theorem}[The kernel $\varphi_p''$ is variation diminishing]\label{thm:kernel-vd}
Let $h\colon(0,\infty)\to\R$ be a finite step function that vanishes for
$r>R$, for some $R>0$, and is not zero almost everywhere. Then
\[
 f(p)=\int_0^\infty\varphi_p''(r)\,h(r)\,dr
\]
has at most $S^-(h)$ zeros in $[0,\infty)$, counted with multiplicity.
\end{theorem}

\noindent The restriction to $p\ge0$ is necessary. For $p<0$ the
function $\varphi_p'=\tanh r\cosh^pr$ increases and then decreases, so
$\varphi_p'(a)=\varphi_p'(b)$ for suitable $0<a<b$, and $h=\mathbf
1_{(a,b)}$ has $S^-(h)=0$, while
$\int_0^\infty\varphi_p''h\,dr=\varphi_p'(b)-\varphi_p'(a)=0$.

\begin{proposition}[Tail bound]\label{prop:hyp-tail}
Suppose the center map is nonconstant. Let $P\in\Hh^2$, and let
$\boldsymbol e\in T_P\Hh^2$ be a unit vector with
$N_{\boldsymbol e}\not\equiv0$. Then
\[
 \kappa(P)\le S^-(N_{\boldsymbol e})\le V_{\boldsymbol e}(P).
\]
\end{proposition}

\begin{proof}
By Lemma~\ref{lem:order-contact}, at each finite parameter $p_0$ with
$C_{p_0}=P$, the projected force $\langle F_p(P),\boldsymbol e\rangle$
vanishes to order at least $\operatorname{ord}_{p_0}(P)$. By the tail
identity~\eqref{eq:hyp-tail}, this projection is the transform of
$N_{\boldsymbol e}$. Since $N_{\boldsymbol e}\not\equiv0$,
Theorem~\ref{thm:kernel-vd} gives $\kappa(P)\le S^-(N_{\boldsymbol e})$. 

For the second inequality, read the values of $N_{\boldsymbol e}$ from
the outermost shell inward. They are the partial sums, starting at $0$,
of the numbers $\langle\boldsymbol S_j(P),\boldsymbol e\rangle$ taken
in the order $j=s,s-1,\ldots,1$. Among the terms forming the first
nonzero partial sum, at least one has its sign. Each time the nonzero
partial sums change sign, at least one of the terms added since the
last partial sum of the old sign has the new sign: their sum is the
difference of the two partial sums. The terms so obtained occur in
order and alternate in sign. Hence the partial sums change sign no
more often than the terms themselves. Reversing the order of the
terms does not change their number of sign changes, so
$S^-(N_{\boldsymbol e})\le V_{\boldsymbol e}(P)$.
\end{proof}

\noindent For points on a geodesic through $P$, with unit tangent
$\boldsymbol u$ at $P$, every direction $\vp_i$ is $\pm\boldsymbol u$, so
the tail path moves along the line $\R\boldsymbol u$ by integer steps and
$N_{\boldsymbol u}$ is an integer walk. This already shows how much
smaller the tail bound can be than the sign bound. If no point is at $P$
and each shell contains one point, on alternating sides of $P$ in order
of distance, as in the proof of Theorem~\ref{thm:n-minus-two}, then the
shell vectors $\pm\boldsymbol u$ alternate and $V_{\boldsymbol u}(P)=n-1$,
the largest possible value. If $\boldsymbol u$ points toward the farthest
point, however, $N_{\boldsymbol u}$ takes only the values $0$ and $1$, so
$S^-(N_{\boldsymbol u})=0$, and Proposition~\ref{prop:hyp-tail} shows
again that $P$ is never the center, as did the pairing argument in the proof of Theorem~\ref{thm:n-minus-two} (compare Figure~\ref{fig:tail}, where the same
happens for a tail path that is not contained in a line). In general,
leaving $0$ costs the walk at least one unit of variation and every
subsequent change of sign at least two, which gives the collinear bound
of the next theorem.

\begin{theorem}[Points on a geodesic]\label{thm:collinear-multiplicity}
Let $\mathcal A$ consist of $n\ge3$ points on a geodesic $\gamma$, counted
with multiplicity, and suppose that the center map is
nonconstant. Then $C_p\in\gamma$ for every $p\in[0,\infty]$, and for every
$P\in\Hh^2$
\[
 \kappa(P)+\mathbf 1_{\{C_\infty=P\}}\le\Bigl\lfloor\frac{n-1}2\Bigr\rfloor .
\]
\end{theorem}

\begin{proof}
The geodesic is closed and convex, so $C_p\in\gamma$ by
Proposition~\ref{prop:convex-hull}. For $P\notin\gamma$ the left side is
zero, so we may assume $P\in\gamma$. Let
$\boldsymbol u$ be a unit tangent vector of $\gamma$ at $P$ and $n'$ the
number of points different from $P$. Each of them has
$\vp_i=\pm\boldsymbol u$, so $N_{\boldsymbol u}$ is an integer
step function that vanishes beyond the largest distance $R$, with inward jump
$\langle\boldsymbol S_j(P),\boldsymbol u\rangle$ at $r_j$, and with total
variation $\sum_j|\langle\boldsymbol S_j(P),\boldsymbol u\rangle|\leq n'$. Since the center map is nonconstant, $P$ is not a
universal center. It follows that some shell vector is a nonzero multiple of
$\boldsymbol u$, so $N_{\boldsymbol u}\not\equiv0$. Assume
$N_{\boldsymbol u}$ changes sign $V$ times and let
$\rho_0>\rho_1>\cdots>\rho_V$ be the radii at which it takes nonzero integer values
of alternating signs. The total variation satisfies
$n'\geq \sum_j|\langle\boldsymbol S_j(P),\boldsymbol u\rangle| \geq |N_{\boldsymbol u}(\rho_0)|+\sum_l|N_{\boldsymbol u}(\rho_l)-
N_{\boldsymbol u}(\rho_{l-1})|\ge1+2V$. Hence $V\le\lfloor(n'-1)/2\rfloor$,
and Proposition~\ref{prop:hyp-tail} gives the bound. If $C_\infty=P$, then Proposition~\ref{prop:equilibrium} says that $0$ lies in the convex hull of the directions $\pm\boldsymbol u$ to the farthest points, so both directions must occur. Let $a\ge1$
be the number of points in the outermost shell in the direction $\boldsymbol u$ and $c\ge1$ in the direction
$-\boldsymbol u$. The jump $a-c$ of this shell satisfies $|a-c|\le a+c-2$. This reduces the total variation to at most $n'-2$ and the rest is as before: $n'-2\geq 1+2V$, so
$\kappa(P)\le V \le\lfloor(n'-3)/2\rfloor$ and
$\kappa(P)+1\le\lfloor(n'-1)/2\rfloor\le\lfloor(n-1)/2\rfloor$.
\end{proof}

\begin{remark}[Higher dimensions]\label{rem:scope}
The results of this section extend to $\Hh^d$ for every $d\ge2$. The
geometric preliminaries, analyticity of the center map, and
the multiplicity argument remain valid in every dimension. The
shell-span bound still reads
\[
 \kappa(P)\le m-k,
\]
where $m$ is the number of nonzero shell vectors and $k$ is the dimension
of their span in $T_P\Hh^d$, rather than the dimension of the geodesic
hull of the configuration, since shells may cancel. Choose a unit vector
$\boldsymbol e$ orthogonal to $k-1$ linearly independent shell vectors but
not orthogonal to the span of all $m$ unbalanced shell vectors. The
projected sequence $\langle\boldsymbol S_j,\boldsymbol e\rangle$ then has at
most $m-k+1$ nonzero terms, so Theorem~\ref{thm:sign-bound}(i) gives
$\kappa(P)\le m-k$. For $k\ge3$, this improves
the estimate $m-2$ available when one only uses $k\ge2$.
The proof of Theorem~\ref{thm:n-minus-two} is unchanged: the exceptional
case $k=1$, $m=n$ again reduces to points on a geodesic. The tail
identity holds in every tangent space, and the scalar kernel
$\varphi_p''$ is unchanged, so the tail bound and the collinear bound
also extend. The binary-form applications in Section~\ref{sec:forms} are
specific to $d=2$.
\end{remark}

\subsection{\texorpdfstring{Proof of Theorem~\ref{thm:kernel-vd}}{Proof of the variation-diminishing theorem}}\label{subsec:proof-vd}

The proof rests on two lemmas. The first bounds the zeros of linear
combinations of the radial weights $\varphi_p'=w_p$ and their
derivatives in $p$; the second is P\'olya's Wronskian criterion.

\begin{lemma}\label{lem:dirichlet}
Let $0\le p_1<\cdots<p_t$ and let $d_1,\ldots,d_t\ge0$ be integers, and
let $J=\sum_{\ell}(d_{\ell}+1)$. The $J+1$ functions
\[
 1\quad\text{and}\quad
 \partial_p^k\varphi_p'(r)\big|_{p=p_{\ell}}\qquad(1\le \ell\le t,\ 0\le k\le d_{\ell})
\]
are linearly independent on $(0,\infty)$, and every nontrivial linear
combination of them has at most $J$ zeros in $(0,\infty)$, counted with
multiplicity.
\end{lemma}

\begin{proof}
Let $\lambda=\log\cosh r$ and write
$v(\lambda)=\tanh r=(1-e^{-2\lambda})^{1/2}$.
In this variable the radial weight is
$\varphi_p'(r)=v(\lambda)e^{p\lambda}$, so each differentiation
in $p$ contributes a factor of $\lambda$:
\[
 \left.\partial_p^k\varphi_p'(r)\right|_{p=p_\ell}
 =\lambda^k v(\lambda)e^{p_\ell\lambda}.
\]

Thus, after division by $v(\lambda)$, a linear combination of the
functions in the statement takes the form
\[
 f(\lambda)=\frac{c_0}{v(\lambda)}
       +\sum_{\ell=1}^t P_\ell(\lambda)e^{p_\ell\lambda},
 \qquad \deg P_\ell\le d_\ell.
\]
The change of variable is an analytic diffeomorphism of $(0,\infty)$
onto itself, and $v$ is positive there, so neither operation changes
the zeros or their multiplicities.

Write $D=d/d\lambda$. For any real $a$,
\[
 (D-a)g=e^{a\lambda}\bigl(e^{-a\lambda}g\bigr)'.
\]
It follows that if $g$ has $N\ge1$ zeros counted with multiplicity,
then $(D-a)g$ has at least $N-1$: differentiation retains all but
one order at each zero, and Rolle's theorem gives a further zero
between consecutive distinct zeros. In particular, if applying $K$
such factors to $g$ gives a nowhere-zero function, then $g$ has
at most $K$ zeros counted with multiplicity.

Suppose first that $c_0\ne0$. Consider the operator
\[
 L=\prod_{\ell=1}^t(D-p_\ell)^{d_\ell+1},
\]
which has order $J$. Since
$(D-p_\ell)(P e^{p_\ell\lambda})=P'e^{p_\ell\lambda}$,
$L$ annihilates every polynomial-exponential term in $f$.
To evaluate the remaining term, use the binomial expansion
\[
 \frac1{v(\lambda)}
 =\sum_{m=0}^\infty a_m e^{-2m\lambda},
 \qquad a_m=4^{-m}\binom{2m}{m}>0.
\]
On every interval $\lambda\ge\delta>0$ the series may be
differentiated term by term, since the exponential decay dominates
every power of $m$. Hence
\[
 Lf(\lambda)=(-1)^Jc_0
 \sum_{m=0}^\infty a_m
 \prod_{\ell=1}^t(2m+p_\ell)^{d_\ell+1}e^{-2m\lambda}.
\]
All terms in the sum are nonnegative because $p_\ell\ge0$, and
every term with $m\ge1$ is positive. Hence $Lf$ is nowhere zero,
so $f$ has at most $J$ zeros.

Suppose now that $c_0=0$ and that some $P_j$ is nonzero.
Apply the factors $(D-p_\ell)^{d_\ell+1}$ for all $\ell\ne j$.
They annihilate the terms with $\ell\ne j$ and leave a term
$Q(\lambda)e^{p_j\lambda}$ with $\deg Q=\deg P_j$. Indeed,
\[
 (D-p_\ell)(P e^{p_j\lambda})
 =\bigl(P'+(p_j-p_\ell)P\bigr)e^{p_j\lambda},
\]
and $p_j-p_\ell\ne0$, so the polynomial degree is preserved.
A further $\deg P_j$ applications of $D-p_j$ leave a nonzero
constant multiple of $e^{p_j\lambda}$. The total number of factors in the differential operator applied is
\[
 \sum_{\ell\ne j}(d_\ell+1)+\deg P_j\le J-1,
\]
so in this case $f$ has at most $J-1$ zeros.

In both cases a combination with coefficients not all zero cannot
vanish identically, which proves linear independence. The two cases
also give the asserted bound of $J$ zeros, counted with multiplicity.
\end{proof}

\begin{lemma}[P\'olya's Wronskian criterion]\label{lem:polya}
Let $u_1,\ldots,u_M$ be smooth real-valued functions on an open
interval $I$. If $W(u_1,\ldots,u_j)$ is nowhere zero on $I$ for
each $j=1,\ldots,M$, then every nontrivial linear combination of
$u_1,\ldots,u_M$ has at most $M-1$ zeros in $I$, counted with
multiplicity.
\end{lemma}

\noindent See P\'olya~\cite{Polya} and, for this formulation,
Novikov and Yakovenko~\cite[Theorem~5, p.~180]{NovikovYakovenko}.

\begin{proof}[Proof of Theorem~\ref{thm:kernel-vd}]
Let $V=S^-(h)$ and $M=V+1$. Since $h$ has bounded support and
\[
 \varphi_p''(r)=e^{(p-2)\log\cosh r}(1+p\sinh^2r),
\]
the integral defining $f$ is real analytic in $p$ on $\R$ and may
be differentiated under the integral sign. In particular, a zero at
$p=0$ is counted with its usual multiplicity. If the asserted bound
fails, there are parameters $0\le p_1<\cdots<p_t$ and positive
integers $\mu_1,\ldots,\mu_t$ with $\sum_{\ell=1}^t\mu_\ell=M$
such that
\[
 f^{(k)}(p_\ell)=0
 \qquad (1\le\ell\le t,\ 0\le k<\mu_\ell).
\]

For each such pair $(\ell,k)$, let
\[
 U_{\ell,k}(r)=
 \left.\partial_p^k\varphi_p'(r)\right|_{p=p_\ell},
 \qquad
 u_{\ell,k}(r)=U_{\ell,k}'(r)
 =\left.\partial_p^k\varphi_p''(r)\right|_{p=p_\ell},
\]
where the prime denotes differentiation in $r$.
List the $U_{\ell,k}$ in order of increasing $\ell$ and, within each
block, increasing $k$. Write the resulting list as $U_1,\ldots,U_M$
and let $u_i=U_i'$. For $1\le j\le M$, expansion along the column
of the constant function $1$ yields
\[
 W(u_1,\ldots,u_j)=W(1,U_1,\ldots,U_j).
\]
If this determinant vanished at some $r>0$, a nontrivial linear
combination of $1,U_1,\ldots,U_j$ would vanish there together with
its first $j$ derivatives. This contradicts
Lemma~\ref{lem:dirichlet} with $J=j$: in the initial list
$U_1,\ldots,U_j$, each parameter that occurs is represented by
consecutive derivative orders starting at zero. Hence every initial
Wronskian is nowhere zero. By Lemma~\ref{lem:polya}, every
nontrivial linear combination of $u_1,\ldots,u_M$ has at most
$M-1=V$ zeros in $(0,\infty)$, counted with multiplicity.

Choose $\rho_1<\cdots<\rho_V$, one between each successive pair of
nonzero values of $h$ of opposite sign. On the intervals of
$(0,\infty)$ determined by these points, $h$ is alternately
nonnegative and nonpositive almost everywhere. The conditions
$\psi(\rho_i)=0$, $1\le i\le V$, impose $V$ homogeneous linear
equations on the $V+1$ coefficients of
\[
 \psi=\sum_{\ell,k}a_{\ell,k}u_{\ell,k}.
\]
Choose a solution with coefficients not all zero. Since $\psi$ has
at most $V$ zeros counted with multiplicity, the $\rho_i$ are its
only zeros and each is simple. Hence $\psi$ changes sign at each
$\rho_i$ and nowhere else. After replacing $\psi$ by $-\psi$ if
necessary, we have $\psi h\ge0$ almost everywhere. Since $h$ is not
zero almost everywhere and $\psi$ has only finitely many zeros,
\[
 \int_0^\infty \psi(r)h(r)\,dr>0.
\]
On the other hand, differentiation under the integral sign yields
\[
 \int_0^\infty \psi(r)h(r)\,dr
 =\sum_{\ell=1}^t\sum_{k=0}^{\mu_\ell-1}
   a_{\ell,k}f^{(k)}(p_\ell)=0.
\]
This contradiction proves the theorem.
\end{proof}
\section{Sharpness at prescribed parameters}\label{sec:sharp}

In this section we show that the bound $n-2$ of
Theorem~\ref{thm:n-minus-two} and the collinear bound
$\lfloor(n-1)/2\rfloor$ of Theorem~\ref{thm:collinear-multiplicity}
are both sharp. In any neighborhood of a given point $P$, we construct
configurations for which the finite solutions of $C_p=P$ are precisely
the prescribed parameters, with the prescribed multiplicities. The
endpoint may also be prescribed and contributes one to the total, which
equals the corresponding bound.

We use the tail identity in reverse. We first choose a tail count
whose transform has the prescribed finite zeros, then place points at
its jumps to realize it. The choice rests on the Hobby--Rice theorem
\cite{HobbyRice,Pinkus}: for integrable functions $k_1,\ldots,k_m$
on $(0,R)$, there is a step function $\varepsilon$ with values $\pm1$
and at most $m$ sign changes such that
$\int_0^R k_i(r)\varepsilon(r)\,dr=0$ for every $i$.

\begin{lemma}\label{lem:canonical}
Let $R>0$ and $m\ge1$, and let $k_1,\ldots,k_m$ be integrable
functions on $(0,R)$ with the following property: if a finite
step function $h$ on $(0,R)$ satisfies
\[
  \int_0^R k_i(r)h(r)\,dr=0
  \qquad(i=1,\ldots,m),
\]
then $h=0$ almost everywhere or $S^-(h)\ge m$.
Then there exist unique points
\[
  0=t_0<t_1<\cdots<t_m<t_{m+1}=R
\]
such that the step function
\[
  \varepsilon(r)=(-1)^{m-j},
  \qquad t_j<r<t_{j+1},\quad j=0,\ldots,m,
\]
satisfies
\[
  \int_0^R k_i(r)\varepsilon(r)\,dr=0
  \qquad(i=1,\ldots,m).
\]
The values of $\varepsilon$ at the cut points are immaterial.
\end{lemma}

\begin{proof}
The Hobby--Rice theorem gives a step function $\eta$ with values $\pm1$
and at most $m$ sign changes such that $\int_0^Rk_i\eta\,dr=0$ for all
$i$. By the hypothesis, $\eta$ has exactly $m$ sign changes, at points
$0<t_1<\cdots<t_m<R$. Then $\varepsilon=\pm\eta$, with the sign chosen so
that $\varepsilon=1$ on $(t_m,R)$, has the stated form.

For uniqueness, suppose that $\varepsilon$ and $\widetilde\varepsilon$
are two such sign functions, and let
$h=\varepsilon-\widetilde\varepsilon$. Both equal $1$ near $R$ and,
since each changes sign exactly $m$ times, both equal $(-1)^m$ near
$0$. Hence $h$ starts and ends at zero, and its only possible nonzero
values are $2$ and $-2$. Ignoring values at the breakpoints, each sign
function has total variation $2m$, so the total variation of $h$ is
at most $4m$. If $h$ is nonzero on an interval, its total variation
is at least $2+4S^-(h)+2$: leaving zero costs $2$, each change of
sign costs at least $4$, and returning to zero costs $2$. Hence
$S^-(h)\le m-1$. But $\int_0^R k_i h\,dr=0$ for every $i$, so the
hypothesis of the lemma forces $h=0$ almost everywhere. It follows that
$\varepsilon$ and $\widetilde\varepsilon$ agree on every open interval
between their combined cut points. A sign change of one must then be
a sign change of the other, so their cut points coincide.
\end{proof}

\noindent Let $\mathcal Q$ be a finite multiset in $[0,\infty)$, with
$|\mathcal Q|=m$. For each parameter $p$ occurring in $\mathcal Q$ with
multiplicity $\mu$, take the $\mu$ functions
\[
 \partial_s^k\varphi_s''(r)\big|_{s=p},\qquad 0\le k<\mu .
\]
By Theorem~\ref{thm:kernel-vd}, these $m$ functions have the property
required in Lemma~\ref{lem:canonical}: if a finite step function $h$ on
$(0,R)$ is orthogonal to all of them then, differentiating under the
integral sign, $f(p)=\int_0^R\varphi_p''h\,dr$ has a zero of
multiplicity at least $\mu$ at each such $p$, hence at least $m$ zeros
counted with multiplicity, so $h=0$ almost everywhere or $S^-(h)\ge m$.

\begin{theorem}[Prescribed parameters and multiplicities]\label{thm:prescribed}
Let $P\in\Hh^2$, $R>0$ and $n\ge3$, and let $\mathcal P$ be a finite
multiset in $[0,\infty]$ in which $\infty$ occurs at most once. In each
of the following cases there is a configuration of $n$ points, counted
with multiplicity, in the closed ball of radius $R$ about $P$, whose
fiber multiset over $P$ is $\mathcal P$:
\begin{enumerate}[label=\textup{(\roman*)}]
\item $|\mathcal P|=\lfloor(n-1)/2\rfloor$, with all configuration points on a geodesic through $P$;
\item $|\mathcal P|=n-2$.
\end{enumerate}
Thus \textup{(i)} attains the collinear bound of
Theorem~\ref{thm:collinear-multiplicity}, and \textup{(ii)} attains
the general bound of Theorem~\ref{thm:n-minus-two}.
\end{theorem}

\begin{proof}
Let $\delta=1$ if $\infty\in\mathcal P$ and $\delta=0$ otherwise, and
let $\mathcal Q$ be the multiset of finite parameters in $\mathcal P$, so
that $m+\delta=|\mathcal P|$ with $m=|\mathcal Q|$.
Fix orthonormal $\boldsymbol u,\boldsymbol e\in T_P\Hh^2$. If $m\ge1$, let
$\varepsilon$, with sign changes $t_1<\cdots<t_m$, be given by
Lemma~\ref{lem:canonical} for the $m$ functions associated with
$\mathcal Q$ above; if $m=0$, let
$\varepsilon\equiv1$ on $(0,R)$. Let $c=1$ in case (i) and $c=\frac12$ in
case (ii), and let $N=c\,\varepsilon$ on $(0,R)$ and $N=0$ on
$(R,\infty)$, so that $S^-(N)=m$. By the tail identity \eqref{eq:hyp-tail},
$F_p(P)=\int_0^\infty\varphi_p''(r)\,\boldsymbol T_P(r)\,dr$, and
$\boldsymbol T_P=N_{\boldsymbol u}\boldsymbol u+N_{\boldsymbol e}\boldsymbol e$.
Hence a configuration whose tail counts are $N_{\boldsymbol u}=N$ and
$N_{\boldsymbol e}=0$ has net $p$-force $f(p)\boldsymbol u$ at $P$, with
$f(p)=\int_0^R\varphi_p''N\,dr$. This function vanishes at each
$p\in\mathcal Q$ to order at least its multiplicity.
By Theorem~\ref{thm:kernel-vd} it has at most $m$ zeros in $[0,\infty)$,
counted with multiplicity, so it has no other zeros and these orders are
exact. By Proposition~\ref{prop:equilibrium} and
Lemma~\ref{lem:order-contact}, $C_p=P$ for finite $p$ exactly when
$p\in\mathcal Q$, and the multiplicity of each such solution is its
multiplicity in $\mathcal Q$.

It remains to realize $N$. Inside the ball, place at distance $t_j$ in
the direction $(-1)^{m-j+1}\boldsymbol u$, for $j=1,\ldots,m$, a point of
multiplicity $2$ in case (i) and a single point in case (ii); these
produce the jumps of $N$ at the $t_j$. At distance $R$, the unit directions must sum to $c\,\boldsymbol u$ to give the required inward jump at $R$. Their convex hull determines whether $P$ is the circumcenter. We choose the directions as follows:
\[
\begin{array}{c|c|c}
 & \delta=0 & \delta=1\\\hline
\textup{(i)} & \boldsymbol u & \boldsymbol u,\ \boldsymbol u,\ -\boldsymbol u\\
\textup{(ii)} & \frac14\boldsymbol u\pm\frac{\sqrt{15}}4\boldsymbol e
 & \boldsymbol u,\ -\frac14\boldsymbol u\pm\frac{\sqrt{15}}4\boldsymbol e
\end{array}
\]
In case (ii) these are $n$ points; in case (i) they are $n$ points when
$n$ is odd, and we add the point $P$ when $n$ is even. The tail counts
are $N_{\boldsymbol u}=N$ and $N_{\boldsymbol e}=0$, and in case (i) the
points lie on one geodesic. The farthest points are those at distance
$R$. For $\delta=0$ their directions have positive projection on
$\boldsymbol u$, so $0$ is not in their convex hull and $C_\infty\ne P$
by Proposition~\ref{prop:equilibrium}. For $\delta=1$ we have
$\frac12\boldsymbol u+\frac12(-\boldsymbol u)=0$ in case (i) and
$\frac15\boldsymbol u+\frac25\bigl(-\frac14\boldsymbol u
+\frac{\sqrt{15}}4\boldsymbol e\bigr)+\frac25\bigl(-\frac14\boldsymbol u
-\frac{\sqrt{15}}4\boldsymbol e\bigr)=0$ in case (ii), so
$C_\infty=P$.
\end{proof}

\noindent For fixed $R$ and prescribed fiber multiset, the cut radii in this
construction are unique by Lemma~\ref{lem:canonical}; the realizing
configuration need not be unique, even up to isometries fixing $P$.

\begin{remark}[Explicit configurations for odd parameters]\label{rem:odd-explicit}
For the multiset $\mathcal Q=\{1,3,\ldots,2m-1\}$, the cut points
$t_j$ in the proof of Theorem~\ref{thm:prescribed} can be chosen
explicitly as
\[
 \sinh t_j=\sinh R\,\cos\frac{(m+1-j)\pi}{2m+1},\qquad j=1,\ldots,m ,
\]
so the numbers $\sinh t_j/\sinh R$ are the positive zeros of the Chebyshev
polynomial $U_{2m}$. Indeed, with $y=\sinh r$ the radial weights are
$w_{2i-1}(r)=\tanh r\,\cosh^{2i-1}r=y(1+y^2)^{i-1}$, so the forms
$\varphi_{2i-1}''\,dr=dw_{2i-1}$, $1\le i\le m$, span the forms
$y^{2l}\,dy$, $0\le l<m$. Let $y=\sinh R\cos\theta$. For $0\le k<2m$ the
form $y^k\,dy$ is a combination of $\sin(l\theta)\,d\theta$ with
$1\le l\le2m$, and $\operatorname{sgn}\sin((2m+1)\theta)$ is orthogonal to
these on $(0,\pi)$, as in the proof of
Theorem~\ref{thm:Lq-collinear}(ii). Since $2m+1$ is odd, this sign
function is invariant under $\theta\mapsto\pi-\theta$, and hence is even
in $y$, so on
$0<y<\sinh R$ it is orthogonal to the forms $y^{2l}\,dy$; it equals $1$
near $y=\sinh R$ and changes sign exactly at the points above. For example,
for $m=2$ and $\sinh R=4$, the five points on a geodesic through $P$ at
signed distances $\operatorname{arsinh}4$, $-\operatorname{arsinh}(\sqrt5+1)$
and $\operatorname{arsinh}(\sqrt5-1)$, the last two with multiplicity $2$,
have $C_p=P$ exactly for $p\in\{1,3\}$.
\end{remark}

\begin{remark}[Higher dimensions]\label{rem:scope-constructions}
The configurations of Theorem~\ref{thm:prescribed} lie in a totally geodesic plane through $P$. In
$\Hh^d$ the centers of a configuration contained in such a plane lie in
the plane, by Proposition~\ref{prop:convex-hull}, and therefore agree with
the centers computed in the plane, which is isometric to $\Hh^2$. Hence
the construction, and the sharpness of the bounds of
Section~\ref{sec:rigidity}, hold in $\Hh^d$ for every $d\ge2$.
\end{remark}

\section{Applications to binary forms}\label{sec:forms}

For a totally complex real binary form $F$ of degree $2n$, the $p$-zero
$\xi_p(F)$ is the center $C_p$ of its $n$ upper-half-plane roots
(Subsection~\ref{subsec:binary-forms}). Throughout this section roots in
$\Hh^2$ are counted with multiplicity, and the Julia extension of
Subsection~\ref{subsec:binary-forms} is used when all upper-half-plane roots coincide. Conversely, every configuration $\mathcal A$ of $n$ points in $\Hh^2$ is
the multiset of upper-half-plane roots of the totally complex real form
$F(X,Z)=\prod_{\alpha\in\mathcal A}(X-\alpha Z)(X-\overline\alpha Z)$ of
degree $2n$. We translate the bounds of
Section~\ref{sec:rigidity} and the constructions of
Section~\ref{sec:sharp} into coincidence thresholds for $p$-zeros,
classify the forms of low degree with a parameter-independent zero, and describe when a reducing matrix remains valid as \(p\) varies.

\subsection{Coincidence thresholds}\label{subsec:thresholds}

By Theorem~\ref{thm:n-minus-two}, if sufficiently many distinct
parameters give the same $p$-zero of $F$, then $\xi_p(F)$ is
independent of $p$.

\begin{corollary}[Coincidence of $p$-zeros]\label{cor:forms-n-minus-one}
If the $p$-zeros of a totally complex real binary form $F$ of degree $2n$ take
the same value at $n-1$ distinct parameters in $[0,\infty]$, then
$\xi_p(F)$ is the same point for every $p\in[0,\infty]$. For sextics any
two parameters suffice, and for octics any three.
\end{corollary}

\begin{proof}
For $n\le2$ the center map is constant
(Remark~\ref{rem:two-points}). For $n\ge3$, apply
Theorem~\ref{thm:n-minus-two} to the $n$ upper-half-plane roots.
\end{proof}

\noindent For $n\ge3$ the threshold $n-1$ is sharp. Given a point
$P\in\Hh^2$ and any $n-2$ distinct parameters in $[0,\infty]$,
Theorem~\ref{thm:prescribed}(ii) gives a form of degree $2n$ whose
$p$-zero equals $P$ exactly at those parameters, with all its
upper-half-plane roots arbitrarily close to $P$.

In particular, for every $n\ge4$ and every pair of distinct parameters
$p_1,p_2\in[0,\infty]$, there is a form $F$ of degree $2n$ such that
$\xi_{p_1}(F)=\xi_{p_2}(F)$ but $p\mapsto\xi_p(F)$ is nonconstant.
Taking $p_1=0$ and $p_2=1$, we see that degree~$8$ is the first
degree in which Julia's zero and the hyperbolic zero can coincide
without all $p$-zeros coinciding; see Example~\ref{ex:integral-octic}
for an explicit integral octic.

For collinear roots, the sharp threshold is instead
$\lfloor(n-1)/2\rfloor+1$, by Theorems~\ref{thm:collinear-multiplicity}
and~\ref{thm:prescribed}(i).

\subsection{Parameter-independent zeros of sextics and octics}\label{subsec:low-degree}

We now describe the root configurations for which the $p$-zero is
independent of $p$ in degrees $4$, $6$, and $8$. For quartics, every
$p$-zero is the hyperbolic midpoint of the two upper-half-plane roots
(Remark~\ref{rem:two-points}). For the pair of parameters $p=0,1$,
this observation, Theorem~\ref{thm:sextic-universal}, and
Corollary~\ref{prop:octic-universal} recover the coincidence results
of \cite[Propositions~11--13]{EleziEquilibrium} for quartics, sextics,
and collinear octics, respectively. The proofs here are self-contained.

\begin{theorem}[Sextics with coincident $p$-zeros]\label{thm:sextic-universal}
Let $\alpha_1,\alpha_2,\alpha_3\in\Hh^2$ be the upper-half-plane roots
of a totally complex real binary sextic $F$. The following are equivalent:
\begin{enumerate}[label=\textup{(\roman*)}]
\item $\xi_p(F)=\xi_{p'}(F)$ for some distinct $p,p'\in[0,\infty]$;
\item either $\alpha_1,\alpha_2,\alpha_3$ are the vertices of an equilateral
hyperbolic triangle, or they are collinear and one of them is the hyperbolic
midpoint of the other two;
\item $\xi_p(F)$ is independent of $p\in[0,\infty]$.\end{enumerate}
\end{theorem}

\begin{proof}
By Corollary~\ref{cor:forms-n-minus-one}, coincidence at two distinct
parameters suffices for a sextic, so \textup{(i)} implies \textup{(iii)}.

Assume \textup{(iii)}, and let $P$ be the common $p$-zero. By
Theorem~\ref{thm:shell}, every shell about $P$ is balanced and therefore
contains at least two roots. Since there are only three roots, there is
at most one shell. If it contains all three roots, their unit directions
from $P$ sum to zero and hence make pairwise angles of $2\pi/3$. The
roots are equidistant from $P$, so the hyperbolic cosine law yields the
same distance between each pair. They therefore form an equilateral
triangle. If the shell contains two roots, their directions are opposite
and their distances from $P$ are equal, so $P$ is their midpoint. The
third root is $P$ itself. If there is no shell, all three roots equal
$P$, and the midpoint condition again holds. This proves \textup{(ii)}.

Conversely, suppose \textup{(ii)} holds. In the equilateral case, rotation
by $2\pi/3$ about the circumcenter permutes the three vertices, so their
unit directions there sum to zero. In the midpoint case, the two outer
roots form a balanced shell about their midpoint, and the root at the
midpoint contributes no force. If all three roots coincide, there are
no shells. Hence in each case Theorem~\ref{thm:shell} gives
\textup{(iii)}. Finally, \textup{(iii)} immediately implies \textup{(i)}.
\end{proof}

For sextics, coincidence at two distinct parameters forces all $p$-zeros
to coincide. For octics, even coincidence of Julia's zero and the
hyperbolic zero does not suffice (Example~\ref{ex:integral-octic}). The
symmetry principle \cite[Proposition~10]{EleziEquilibrium} explains why
these two zeros agree in each of the configurations below, and the shell
criterion shows that these are exactly the octic configurations with a
parameter-independent zero.

\begin{theorem}[Octics with a parameter-independent zero]
\label{prop:octic-classification}
Let $F$ be a totally complex real binary octic, and let $P\in\Hh^2$.
Then $\xi_p(F)=P$ for every $p\in[0,\infty]$ if and only if the four
upper-half-plane roots of $F$, counted with multiplicity, have one of
the following configurations:
\begin{enumerate}[label=\textup{(\roman*)}]
\item The roots can be partitioned into two pairs, each with
      hyperbolic midpoint $P$.
\item One root is $P$, and the other three are the vertices of an
      equilateral hyperbolic triangle with circumcenter $P$.
\end{enumerate}
In \textup{(i)}, the two pairs need not lie on the same geodesic,
and a pair may consist of two copies of $P$.
\end{theorem}

\begin{proof}
By Theorem~\ref{thm:shell}, the zero is independent of $p$ and equal
to $P$ exactly when every shell about $P$ is balanced. Suppose first
that this holds. A shell cannot contain just one root. If it contains
two, their unit directions are opposite, so $P$ is their midpoint.

If a shell contains three roots, their three unit directions sum to
zero and hence make pairwise angles of $2\pi/3$. Since the roots are
equidistant from $P$, they form an equilateral triangle with
circumcenter $P$. The fourth root must be $P$, since otherwise it
would form a shell by itself. This gives \textup{(ii)}.

It remains to consider a shell containing all four roots. Write their
unit directions as $\vp_1,\ldots,\vp_4$. Since
$\vp_1+\vp_2=-\vp_3-\vp_4$, the pairs
$\{\vp_1,\vp_2\}$ and $\{-\vp_3,-\vp_4\}$ on the unit circle in
$T_P\Hh^2$ have the same midpoint. If this midpoint is zero, then
$\vp_1=-\vp_2$ and $\vp_3=-\vp_4$. If it is nonzero, it determines
the chord uniquely: the chord, possibly degenerate, lies on the line through the midpoint perpendicular to its radius. Hence in either case the four directions
split into two opposite pairs. The roots lie at the same distance
from $P$, so both pairs have hyperbolic midpoint $P$.

In all remaining cases, every shell contains two roots. Pair the roots
on each shell and pair any remaining copies of $P$. This gives
\textup{(i)}.

Conversely, in \textup{(i)} each pair away from $P$ contributes two
opposite directions to the same shell. In \textup{(ii)} the three
directions to the vertices are spaced by $2\pi/3$ and sum to zero.
Hence every shell is balanced in either case, and the conclusion
follows from Theorem~\ref{thm:shell}.
\end{proof}

\begin{example}\label{ex:integral-octic}
The following octic belongs to neither class, although Julia's zero
and the hyperbolic zero coincide:
\[
 F=(X^2+4Z^2)(9X^2+Z^2)
   (3969X^4-7398X^2Z^2+5929Z^4).
\]
Its upper-half-plane roots are
\[
 2i,\qquad i/3,\qquad \frac{\pm5\sqrt{19}+8i}{21}.
\]
At $P=i$, the first two roots lie on the imaginary axis and the last
two are symmetric about it. Hence the net force is vertical. Let
$\boldsymbol e$ be the upward unit vector at $i$. The three shells
have $\cosh r=5/4,\,5/3,\,35/12$, respectively, and their base forces
in the direction $\boldsymbol e$ are $3/5,\,-4/5,\,1/5$. Indeed, the
vertical component of the base force at $i$ from a root $\alpha$ is
$(|\alpha|^2-1)/(|\alpha|^2+1)$. Each root of the symmetric pair has
$|\alpha|^2=11/9$ and so contributes $1/10$. Hence
\[
 \langle F_p(i),\boldsymbol e\rangle
 =\frac35\Bigl(\frac54\Bigr)^p
  -\frac45\Bigl(\frac53\Bigr)^p
  +\frac15\Bigl(\frac{35}{12}\Bigr)^p.
\]
This vanishes at $p=0$ and $p=1$, but equals $5/12$ at $p=2$.
It follows that $\xi_0(F)=\xi_1(F)=i$, whereas $\xi_2(F)\ne i$.
\end{example}

For collinear octics, coincidence at two distinct parameters again
forces the zero to be independent of the parameter.

\begin{corollary}[Collinear octics]\label{prop:octic-universal}
Let $F$ be a totally complex real binary octic whose four
upper-half-plane roots, counted with multiplicity, lie on a hyperbolic
geodesic. For a point $P$ on that geodesic, the following are equivalent:
\begin{enumerate}[label=\textup{(\roman*)}]
\item $\xi_p(F)=\xi_{p'}(F)=P$ for some distinct $p,p'\in[0,\infty]$;
\item the four roots can be paired so that $P$ is the hyperbolic
      midpoint of both pairs;
\item $\xi_p(F)=P$ for every $p\in[0,\infty]$.
\end{enumerate}
\end{corollary}

\begin{proof}
If $p\mapsto\xi_p(F)$ is nonconstant, then
Theorem~\ref{thm:collinear-multiplicity} allows at most
$\lfloor(4-1)/2\rfloor=1$ parameter in $[0,\infty]$ at which
$\xi_p(F)=P$. Hence \textup{(i)} implies \textup{(iii)}, and the
reverse implication is immediate.

By Theorem~\ref{prop:octic-classification}, \textup{(iii)} holds
exactly when the roots can be paired as in \textup{(ii)}, or one
root is $P$ and the other three form an equilateral triangle with
circumcenter $P$. Collinearity rules out the triangle, so
\textup{(iii)} is equivalent to \textup{(ii)}.
\end{proof}

\noindent This rigidity fails in degree $10$:
Theorem~\ref{thm:prescribed}(i), with $n=5$ and the parameters $0$ and
$1$, gives five collinear roots with $C_0=C_1=P$ and a nonconstant
center map.

\subsection{Consequences for reduction}\label{sec:reduction}

Let $\mathcal F$ be the standard closed modular fundamental domain,
\[
 \mathcal F=\{z\in\Hh^2:|\operatorname{Re}z|\le\tfrac12,\ |z|\ge1\}.
\]
For a totally complex real binary form $F$ and $p\in[0,\infty]$, call
$F$ \emph{$p$-reduced} if $\xi_p(F)\in\mathcal F$. A matrix
$g\in\SL_2(\mathbb Z)$ \emph{$p$-reduces} $F$ if $F^g$ is $p$-reduced.
By equivariance, with the substitution convention of
Section~\ref{sec:prelim},
\[
 g\text{ $p$-reduces }F
 \quad\Longleftrightarrow\quad g^{-1}\xi_p(F)\in\mathcal F.
\]
Such a matrix exists because every modular orbit meets $\mathcal F$.
The reducing matrix is unique up to sign when the point
$g^{-1}\xi_p(F)$ lies in the interior of $\mathcal F$. On the boundary,
the side identifications give additional choices.

\begin{corollary}[Parameter-independent reduction]
\label{cor:universal-reduction}
If $P$ is the universal center of the upper-half-plane roots of $F$,
then the set of matrices that $p$-reduce $F$ is the same for every
$p\in[0,\infty]$, namely
\[
 \{g\in\SL_2(\mathbb Z):g^{-1}P\in\mathcal F\}.
\]
In particular, a fixed rule that selects a reducing matrix from the
center alone gives the same transformed form for every $p$.
\end{corollary}

\begin{proof}
Substitute $\xi_p(F)=P$ in the reduction criterion above.
\end{proof}

\noindent Corollary~\ref{cor:universal-reduction} applies to every totally complex
quartic and to the sextics and octics classified in
Theorems~\ref{thm:sextic-universal} and~\ref{prop:octic-classification}.
In degree $2n$, it suffices that $\xi_p(F)$ takes the same value at
$n-1$ distinct parameters in $[0,\infty]$, by
Corollary~\ref{cor:forms-n-minus-one}.

\begin{remark}[Persistence of a reducing matrix]\label{cor:local-reduction}
Suppose $g^{-1}\xi_{p_0}(F)$ lies in the interior of $\mathcal F$ for
some $p_0\in[0,\infty]$, and let $\delta$ be its hyperbolic distance
from the boundary. Since $g^{-1}$ is an isometry, $g$ still $p$-reduces
$F$ whenever $d(\xi_p(F),\xi_{p_0}(F))<\delta$. By
Theorem~\ref{thm:path}, this holds for all $p$ sufficiently close to
$p_0$. In the case $p_0=\infty$, the bound~\eqref{eq:quant-bound} shows
that $g$ $p$-reduces $F$ for every finite
$p>\log n/\log(\cosh\delta)$, where $2n$ is the degree of $F$.

For finite $p$ and $p_0$, the distance between the two zeros is at
most the length of the intervening path, so it suffices that
\[
 \int_{\min\{p,p_0\}}^{\max\{p,p_0\}}
 \Bigl\|\frac{d}{dq}\xi_q(F)\Bigr\|\,dq<\delta.
\]
Here $\|\cdot\|$ is the hyperbolic norm, and the velocity is given by
Corollary~\ref{cor:firstvariation}.
\end{remark}

\noindent In general, different parameters may give different reduced
forms, and neither result compares their heights. For forms with
hyperbolically separated root clusters, this question is taken up in
a sequel in preparation.

\section{Riemannian \texorpdfstring{$L^q$}{Lq} centers}\label{sec:Lq}

Let $M$ be a Hadamard manifold, that is, a complete, simply connected
Riemannian manifold with nonpositive sectional curvature. For $P\ne\alpha$,
write $\vp_\alpha(P)=\exp_P^{-1}(\alpha)/d(P,\alpha)$ for the unit vector
at $P$ pointing toward $\alpha$. Then $M$ is a proper CAT$(0)$ space,
$\exp_P\colon T_PM\to M$ is a diffeomorphism for every $P$, and
$d(\cdot,\alpha)$ is convex, and smooth on $M\setminus\{\alpha\}$ with
gradient $-\vp_\alpha$ \cite[Chapter~II.1--II.4]{BridsonHaefliger}. Examples are
Euclidean spaces, the hyperbolic spaces $\Hh^d$, and symmetric spaces of
noncompact type. 

\noindent For a
multiset $\mathcal A=\{\alpha_1,\ldots,\alpha_n\}$ in $M$ and $1<q<\infty$,
the \emph{Riemannian $L^q$ center}
$E_q=E_q(\mathcal A)$ \cite{Afsari} is the minimizer
of the \emph{$L^q$ potential} \[\Psi_q(P)=\sum_id(P,\alpha_i)^q\] (Lemma~\ref{lem:Lq-basic}); $E_2$ is
the Riemannian center of mass \cite{Karcher}. Equivalently, $E_q$ minimizes the $q$th power mean of the
distances,
\[
 \left(\frac1n\sum_{i=1}^n d(P,\alpha_i)^q\right)^{1/q}.
\]
As $q\to\infty$, this expression converges to
$\max_i d(P,\alpha_i)$, whose minimizer is the circumcenter
$E_\infty$. We call $q\mapsto E_q$, $q\in(1,\infty]$, the
\emph{center path}.

The shell and tail methods of Section~\ref{sec:rigidity} also apply to
Riemannian $L^q$ centers on Hadamard manifolds. We first establish their
basic properties. For $q>1$, define the \emph{net $q$-force} at $P$ by
\begin{equation}\label{eq:Lq-force}
 F_q(P):=\sum_{\alpha_i\ne P}
 d(P,\alpha_i)^{q-1}\,\vp_{\alpha_i}(P).
\end{equation}
Points equal to $P$ contribute zero.

\begin{lemma}\label{lem:Lq-basic}
\begin{enumerate}[label=\textup{(\roman*)}]
\item For $1<q<\infty$, $\Psi_q$ is $C^1$ and strictly convex, with a
unique minimizer $E_q=E_q(\mathcal A)$, and
\begin{equation}\label{eq:Lq-equilibrium}
 \grad\Psi_q(P)=-qF_q(P),
 \qquad
 P=E_q\quad\Longleftrightarrow\quad F_q(P)=0 .
\end{equation}
\item The function $P\mapsto\max_id(P,\alpha_i)$ has a unique minimizer
$E_\infty$. If $\mathcal A$ contains two distinct points, then $P=E_\infty$
if and only if $0$ lies in the convex hull of the vectors
$\vp_{\alpha_i}(P)$ for which $d(P,\alpha_i)$ is maximal.
\item If $\mathcal A$ lies in a closed convex set $K$, then $E_q\in K$ for
every $q\in(1,\infty]$.
\item For every $q>1$,
\[
 d(E_q,E_\infty)\le R\sqrt{n^{2/q}-1},\qquad
 R=\max_id(E_\infty,\alpha_i);
\]
in particular, $E_q\to E_\infty$ as $q\to\infty$.
\end{enumerate}
\end{lemma}

\begin{proof}
(i) Since $q>1$, $d(\cdot,\alpha)^q$ is $C^1$ with gradient
$-q\,d(\cdot,\alpha)^{q-1}\vp_\alpha$, extended by zero at $\alpha$. Indeed, by the
Cartan--Hadamard theorem $d(\cdot,\alpha)$ is smooth on $M\setminus\{\alpha\}$,
and, for $P\ne\alpha$ and $V\in T_PM$, the first variation formula for
distance yields
\[
 \left.\frac{d}{dt}\right|_{t=0}d(\exp_P(tV),\alpha)
 =-\langle V,\vp_\alpha(P)\rangle.
\]
Hence $\grad d(\cdot,\alpha)(P)=-\vp_\alpha(P)$, and the chain rule yields
the stated gradient of $d(\cdot,\alpha)^q$. At $\alpha$ we have $d^q=o(d)$, hence $d^q$ is differentiable at $\alpha$ with zero gradient, and $q\,d^{q-1}\to0$, so the gradient is continuous at $\alpha$. The function $d(\cdot,\alpha)^q$ is
strictly convex along geodesics. To see this, let $x\ne y$ and let
$\gamma\colon[0,1]\to M$ be the geodesic from $x$ to $y$. For $t\in(0,1)$ let
$c_t=\gamma(t)$. Denote $a=d(x,\alpha)$ and $b=d(y,\alpha)$. Then
\[
 d(c_t,\alpha)^q\le\bigl((1-t)a+tb\bigr)^q\le(1-t)a^q+tb^q,
\]
by the convexity of $d(\cdot,\alpha)$ and of $s\mapsto s^q$. If $a\ne b$,
the second inequality is strict, since $s\mapsto s^q$ is strictly convex
for $q>1$. If $a=b$, the first is strict, since the CAT$(0)$ strong convexity inequality for squared distance
\begin{equation}\label{eq:cat0}
 d(c_t,\alpha)^2\le(1-t)a^2+tb^2-t(1-t)\,d(x,y)^2
\end{equation}
\cite[Chapter~II.1]{BridsonHaefliger}, which for $t=\frac12$ is the CN
inequality of Bruhat and Tits, gives $d(c_t,\alpha)<a$. Hence each term
of $\Psi_q$, and so $\Psi_q$ itself, is strictly convex. Since $M$ is proper and $\Psi_q$ is continuous and coercive, a
minimizer exists. Strict geodesic convexity makes it unique. It is
also the only critical point, since every critical point of a $C^1$
geodesically convex function is a minimizer.

(ii) Existence and uniqueness of the circumcenter hold in every complete
CAT$(0)$ space \cite[Proposition~II.2.7]{BridsonHaefliger}. The criterion
is proved as in Proposition~\ref{prop:equilibrium}, with the hyperbolic
cosine law replaced by the CAT$(0)$ angle-comparison inequality
$c^2\ge a^2+b^2-2ab\cos\theta$ \cite[Chapter~II.1]{BridsonHaefliger}. If $0$ lies in
the convex hull of the active directions, those to the farthest points,
and $\boldsymbol u$ is a unit vector, then
$\langle\boldsymbol u,\vp_{\alpha_i}(P)\rangle\le0$ for some active $i$. We obtain
$d(\exp_P(t\boldsymbol u),\alpha_i)^2\ge d(P,\alpha_i)^2+t^2$ for $t\ge0$.
It follows that $r(\exp_P(t\boldsymbol u))>r(P)$ for $t>0$, where
$r=\max_id(\cdot,\alpha_i)$, and $P$ is the minimizer of $r$.
Conversely, suppose that $0$ does not belong to the convex hull $C$
of the active directions. Let $\boldsymbol u\in C$ have smallest
norm. Then $\boldsymbol u\ne0$, and minimality gives
$\langle\boldsymbol u,\boldsymbol v-\boldsymbol u\rangle\ge0$ for every $\boldsymbol v\in C$.
Hence
$\langle\boldsymbol u,\vp_{\alpha_i}(P)\rangle
\ge\|\boldsymbol u\|^2>0$ for every active $i$. Let
$Q_s=\exp_P(s\boldsymbol u/\|\boldsymbol u\|)$. By the first variation formula,
\[
 \frac{d}{ds}\,d(Q_s,\alpha_i)\Big|_{s=0}
 =-\frac{\langle\boldsymbol u,\vp_{\alpha_i}(P)\rangle}{\|\boldsymbol u\|}<0
\]
for every active $i$, while the inactive distances stay below $r(P)$ for small
$s$ by continuity. Hence $r(Q_s)<r(P)$ for small $s>0$, and $P$ is not the
minimizer of $r$.

(iii) The nearest-point projection $\pi$ onto $K$ is $1$-Lipschitz and
fixes $K$ \cite[Proposition~II.2.4]{BridsonHaefliger}, so
$d(\pi P,\alpha_i)\le d(P,\alpha_i)$ for every $i$. Hence $\pi E_q$ is also
a minimizer, and $\pi E_q=E_q$ by uniqueness.

(iv) Let $x=E_q$, $c=E_\infty$ and $r(Q)=\max_id(Q,\alpha_i)$. Then
$r(x)\le\Psi_q(x)^{1/q}\le\Psi_q(c)^{1/q}\le n^{1/q}R$. Let
$\gamma\colon[0,1]\to M$ be the geodesic from $c$ to $x$ and
$c_t=\gamma(t)$ for $t\in(0,1)$. The CAT$(0)$ inequality \eqref{eq:cat0},
applied to this geodesic and to each $\alpha_i$, yields
\[
 R^2\le r(c_t)^2\le(1-t)R^2+t\,r(x)^2-t(1-t)\,d(x,c)^2.
\]
Dividing by $t$ and letting $t\downarrow0$ yields
$d(x,c)^2\le r(x)^2-R^2\le(n^{2/q}-1)R^2$.
\end{proof}

For one point the center path is constant at that point; for two
points it is constant at their midpoint. Hence a nonconstant center
path requires $n\ge3$.

We call $P$ a \emph{universal $L^q$ center} if $E_q=P$ for every
$q\in(1,\infty]$, and let $\nu(P)=\#\{q\in(1,\infty):E_q=P\}$. Unlike $\kappa(P)$, this counts distinct parameters without
multiplicities; see Remark~\ref{rem:Lq-contact}. Shells,
shell vectors and the numbers $V_{\boldsymbol e}(P)$, $m$, $k$ are
defined as in Section~\ref{sec:rigidity}, now with $1\le k\le\dim M$ when
$P$ is not a universal $L^q$ center. With shell radii $r_1<\cdots<r_s$ and
$\lambda_j=\log r_j$, grouping \eqref{eq:Lq-force} by shells gives the
shell representation
\begin{equation}\label{eq:Lq-shell}
 F_q(P)=\sum_{j=1}^sr_j^{q-1}\,\boldsymbol S_j(P)
 =\sum_{j=1}^se^{(q-1)\lambda_j}\,\boldsymbol S_j(P).
\end{equation}

\begin{theorem}[Rigidity of $L^q$ centers]\label{thm:Lq-rigidity}
Let $\mathcal A$ consist of $n$ points of a Hadamard manifold $M$, let
$P\in M$, and let $r_1<\cdots<r_s$ be the positive distances from $P$ to
the points of $\mathcal A$.
\begin{enumerate}[label=\textup{(\roman*)}]
\item $P$ is a universal $L^q$ center if and only if every shell vector
vanishes, and this already follows if $E_q=P$ for $s$ distinct
$q\in(1,\infty)$.
\item If the center path is nonconstant, then
$\nu(P)\le V_{\boldsymbol e}(P)$ whenever the numbers
$\langle\boldsymbol S_j(P),\boldsymbol e\rangle$ are not all zero,
and $\nu(P)\le m-k$.
\item If $n\ge3$ and the center path is nonconstant, then
$\nu(P)+\mathbf 1_{\{E_\infty=P\}}\le n-2$.
\end{enumerate}
\end{theorem}

\begin{proof}
By Lemma~\ref{lem:Lq-basic}(i) and \eqref{eq:Lq-shell}, $E_q=P$ if and
only if $\sum_je^{(q-1)\lambda_j}\boldsymbol S_j(P)=0$.
This is \eqref{eq:shell-exp} with $\tanh r_j$ replaced by $1$ and $p$
replaced by $q-1>0$, and with $\lambda_j=\log r_j$ in place of
$\lambda_j=\log\cosh r_j$. We now follow the proofs of Theorems~\ref{thm:shell},
\ref{thm:sign-bound} and~\ref{thm:n-minus-two}, counting distinct
parameters rather than multiplicities. The required ingredients are
this exponential-sum representation, Lemma~\ref{lem:exponential-zeros},
the circumcenter criterion of Lemma~\ref{lem:Lq-basic}(ii), and the fact
that $r_j^{q-1}$ increases with $j$. One adjustment is needed when
$\dim M>2$: for the shell-span bound in (ii), use the projection
argument of Remark~\ref{rem:scope} instead of the two-dimensional case
distinction in the proof of Theorem~\ref{thm:sign-bound}(ii).

\end{proof}

\begin{remark}[Multiplicities]\label{rem:Lq-contact}
Suppose that the center path is nonconstant, let $P\notin\mathcal A$,
and let $q_0\in(1,\infty)$ be any parameter with $E_{q_0}=P$. Fix an orthonormal basis of $T_PM$ and let $x(q)\in\R^{\dim M}$ be the
coordinates of $\exp_P^{-1}(E_q)$ in this basis, the normal coordinates of
$E_q$ centered at $P$. The function $\Psi_q$ is smooth near $P$, jointly
in $q$, and for $d_i=d(\cdot,\alpha_i)>0$, the Hessian comparison theorem
yields
\[
 \operatorname{Hess}(d_i^q)(V,V)
 \ge q\,d_i^{q-2}\min\{1,q-1\}\,\|V\|^2.
\]
Hence the Hessian of $\Psi_q$ is positive definite away from $\mathcal A$,
and the implicit function theorem makes $q\mapsto E_q$, and hence $x(q)$,
smooth near $q_0\in(1,\infty)$. We define the
\emph{multiplicity} $\operatorname{ord}_{q_0}(P)$ of $q_0$ as a solution of
$E_q=P$ to be the least integer $\ell\ge1$ for
which $x^{(\ell)}(q_0)\ne0$, or $\infty$ if all derivatives
vanish; it does not depend on the
basis. We claim that it is finite, and that the bounds of
Theorem~\ref{thm:Lq-rigidity}(ii)--(iii) remain valid with $\nu(P)$
replaced by $\sum_{1<q<\infty,\ E_q=P}\operatorname{ord}_q(P)$, the endpoint
still counted once. Indeed, the argument of
Lemma~\ref{lem:order-contact} gives
\[
 qF_q(P)=K(q)x(q),
\]
where $K(q)$ is smooth and invertible near $q_0$.
Since $F_q(P)$ is real analytic in $q$ and not identically zero,
it vanishes at $q_0$ to a finite order, and
$x(q)=qK(q)^{-1}F_q(P)$ vanishes to the same order.
Hence $\operatorname{ord}_{q_0}(P)$ is finite and equals the order
of vanishing of $F_q(P)$. Every nonzero scalar projection of
$F_q(P)$ vanishes to at least this order, so
Lemma~\ref{lem:exponential-zeros} gives the same bounds with
multiplicities included. The same argument, using
Lemma~\ref{lem:laplace-zeros}, shows that the collinear bound in
Theorem~\ref{thm:Lq-collinear}(i) below also holds with multiplicities
whenever $P\notin\mathcal A$.

At data points this conclusion can fail. The hyperbolic potential
$\varphi_p(r)$ is analytic in $r^2$ near $r=0$, with leading term
$r^2/2$, whereas the $L^q$ potential $r^q$ is not twice differentiable
there when $1<q<2$. Hence the finite positive-definite Hessian used in Lemma~\ref{lem:order-contact} is unavailable at such data points, as the following example shows. 
 \end{remark}

\begin{example*}
On the real line, consider the multiset
$\{-1,-1,0,8\}$. Denote its $L^q$ center by $x(q)$. It lies in $[-1,8]$ for all $q>1$. Let $q_0=4/3$ and $A=2\log 8$. By
\eqref{eq:Lq-force}, the signed net force at $x\in[-1,8]$ is
\[
 F_q(x)=(8-x)^{q-1}-2(1+x)^{q-1}-\operatorname{sgn}(x)|x|^{q-1}.
\]
At $x=0$,
\[
 F_q(0)=8^{q-1}-2=A\,(q-q_0)+O((q-q_0)^2),
\]
so $x(q)=0$ exactly for $q=q_0$, and the force has a simple zero
there. Joint continuity of the potential
in $(q,x)$ and uniqueness of its minimizer imply $x(q)\to0$ as
$q\to q_0$.

To determine the displacement, we compare the force at parameters
$q$ and $q_0$, keeping the spatial point fixed. At $q_0$,
\[
 F_{q_0}(x)=-\sqrt[3]{x}+O(|x|).
\]
Differentiating with respect to the parameter yields
\[
 \partial_sF_s(x)
 =(8-x)^{s-1}\log(8-x)
 -2(1+x)^{s-1}\log(1+x)
 -\operatorname{sgn}(x)|x|^{s-1}\log|x|,
\]
where the last term is extended by zero at $x=0$.
As $(s,x)\to(q_0,0)$, the first two terms tend to
$2\log8$ and $0$, respectively. The last also tends to zero:
for $s$ near $q_0$, choose $a>0$ with $s-1\ge a$; then
\[
 |x|^{s-1}|\log|x||
 \le |x|^a|\log|x||\longrightarrow0.
\]
Hence $\partial_sF_s(x)\to A$ as $(s,x)\to (q_0,0)$. Applying the mean
value theorem in the parameter, with $x=x(q)$ held fixed, yields
\[
 0=F_q(x)=F_{q_0}(x)+(q-q_0)\,\partial_sF_s(x)\big|_{s=\sigma_q},
\]
where $\sigma_q$ lies between $q_0$ and $q$. Hence, since
$|x|=o(|x|^{1/3})$,
\[
 q-q_0
 =\frac{\sqrt[3]{x}+O(|x|)}{A+o(1)}
 =\frac{\sqrt[3]{x}}{A}\,(1+o(1)).
\]
Cubing and rearranging yield
\[
 x(q)=(2\log8)^3(q-q_0)^3+o(|q-q_0|^3).
\]
Thus the force has a simple zero, but the displacement has cubic
leading order. Counting $q_0$ with multiplicity $3$ would exceed the
bound $n-2=2$ of Theorem~\ref{thm:Lq-rigidity}(iii).

More generally, for an odd integer $k\ge3$ let
$\mathcal A_k=\{-1,-1,0,2^k\}$ and $q_0=1+1/k$. By the same argument,
\[
 x(q)=(2k\log2)^k(q-q_0)^k+o(|q-q_0|^k).
\]
Thus, with only four data points, the displacement can have
arbitrarily high odd leading order while the force at $0$
has a simple zero.
\end{example*}

\noindent With the tail count $N_{\boldsymbol e}$ of
Subsection~\ref{subsec:tail}, inserting
$d^{q-1}=(q-1)\int_0^d r^{q-2}\,dr$ yields
\begin{equation}\label{eq:tail-identity}
\begin{aligned}
 \langle F_q(P),\boldsymbol e\rangle
 &=(q-1)\int_0^\infty r^{q-2}N_{\boldsymbol e}(r)\,dr\\
 &=(q-1)\int_{-\infty}^\infty
 e^{(q-1)\lambda}N_{\boldsymbol e}(e^\lambda)\,d\lambda,
\end{aligned}
\end{equation}
where the second equality uses $r=e^\lambda$.
Let $h(\lambda)=N_{\boldsymbol e}(e^\lambda)$ and
\[
 L(s)=\int_{-\infty}^{\infty}e^{s\lambda}h(\lambda)\,d\lambda,
 \qquad s>0.
\]
Then $\langle F_q(P),\boldsymbol e\rangle=(q-1)L(q-1)$.
The function $h$ is bounded and
vanishes for $\lambda>\log\max_id(P,\alpha_i)$, but it need not vanish as
$\lambda\to-\infty$, where it tends to $N_{\boldsymbol e}(0^+)$. The integral converges absolutely
for $s=q-1>0$. Provided $h$ is not zero almost everywhere, the next lemma
bounds the zeros of $L$, counted with multiplicity, by the
sign changes of $h$. The same bound applies to the projected
force, since the factor $q-1$ is positive. The induction below multiplies by linear factors and therefore
does not preserve the class of bounded step functions. We allow
polynomial growth and extend the definition of sign changes
beyond step functions.

For a measurable function $h$ on an interval, let
$S^-(h)$ be the least $V$ for which the interval can be cut at $V$ points
into subintervals on which $h$ is alternately $\ge0$ and $\le0$ almost
everywhere, and $S^-(h)=\infty$ if no finite such partition exists. For
finite step functions this agrees with the count of
Subsection~\ref{subsec:tail}.

\begin{lemma}\label{lem:laplace-zeros}
Let $h$ be measurable on $\R$, zero on $(\Lambda,\infty)$ for some
$\Lambda$, with $|h(\lambda)|\le C(1+|\lambda|)^k$, not zero almost
everywhere, and with $S^-(h)<\infty$. Then $L(s)=\int_\R e^{s\lambda}h(\lambda)\,d\lambda$ (i) is smooth on $(0,\infty)$
and (ii) it has at most $S^-(h)$ zeros there, counted with multiplicity.
\end{lemma}

\begin{proof}
(i) Since $h$ vanishes on $(\Lambda,\infty)$ and grows at most
polynomially as $\lambda\to-\infty$, exponential decay makes
$\lambda^j e^{s\lambda}h(\lambda)$ absolutely integrable for
every $s>0$ and integer $j\ge0$, locally uniformly in $s$.
Hence $L$ is smooth, with
\[
 L^{(j)}(s)=\int_{\mathbb R}
 \lambda^j e^{s\lambda}h(\lambda)\,d\lambda.
\]
(ii) We induct on $V=S^-(h)$. If $V=0$, $h$ has one sign and $L$ has no
zeros. If $V\ge1$, let $\mu$ be the first cut point. The function $e^{-s\mu}L(s)$ has the same zeros as $L(s)$ and $\frac{d}{ds}\bigl(e^{-s\mu}L(s)\bigr)=e^{-s\mu}\int_\R
e^{s\lambda}h_1(\lambda)\,d\lambda$ where 
$h_1(\lambda)=(\lambda-\mu)h(\lambda)$ satisfies the hypotheses of the Lemma with
$S^-(h_1)\le V-1$. Apply Rolle's theorem to finish the proof
as in Lemma~\ref{lem:exponential-zeros}.
\end{proof}

\noindent Applied to \eqref{eq:tail-identity},
Lemma~\ref{lem:laplace-zeros} gives
$\nu(P)\le S^-(N_{\boldsymbol e})$ whenever
$N_{\boldsymbol e}\not\equiv0$.
For the sharpness constructions below, fix distinct parameters
$q_1,\ldots,q_m>1$. If a step function $h$ on $(0,R)$, not zero
almost everywhere, is orthogonal to all $m$ kernels $r^{q_i-2}$,
then $L(s)=\int_0^R r^{s-1}h(r)\,dr$ vanishes at the $m$ distinct
points $s=q_i-1$. After the substitution $r=e^\lambda$,
Lemma~\ref{lem:laplace-zeros} gives $S^-(h)\ge m$, as required
by Lemma~\ref{lem:canonical}.

\begin{corollary}[Sharpness of the $L^q$ bound]\label{thm:Lq-sharp}
Let $\dim M\ge2$, $n\ge3$, $P\in M$ and $R>0$. For all
$1<q_1<\cdots<q_{n-2}<\infty$ there are $n$ points of $M$ within distance
$R$ of $P$ such that $E_q=P$ exactly for $q\in\{q_1,\ldots,q_{n-2}\}$; in
particular, the bound of Theorem~\ref{thm:Lq-rigidity}(iii) is attained.
\end{corollary}

\begin{proof}
This is the construction in the proof of Theorem~\ref{thm:prescribed}(ii),
for $\mathcal P$ consisting of distinct finite parameters, with $k_i=r^{q_i-2}$, \eqref{eq:tail-identity} in place of \eqref{eq:hyp-tail},
Lemma~\ref{lem:laplace-zeros} in place of Theorem~\ref{thm:kernel-vd},
and orthonormal $\boldsymbol u,\boldsymbol e\in T_PM$. Each point is
placed as $\alpha=\exp_P(r\boldsymbol v)$, with the radius $r$ and unit
vector $\boldsymbol v$ of that construction, so that $d(P,\alpha)=r$ and
$\vp_\alpha(P)=\boldsymbol v$. Hence $\boldsymbol T_P=\tfrac12\varepsilon\,\boldsymbol u$ on $(0,R)$ and
$\boldsymbol T_P=0$ beyond $R$, where the cut radii of $\varepsilon$ are
now determined by the kernels $r^{q_i-2}$. No totally geodesic surface
is needed. The directions to the farthest points have
positive projection on $\boldsymbol u$, so $E_\infty\ne P$ by
Lemma~\ref{lem:Lq-basic}(ii).
\end{proof}

\noindent The construction in the proof of Theorem~\ref{thm:prescribed}(ii),
for $\mathcal P$ containing $\infty$, carries over in the same way: one of the $n-2$ parameters may
be replaced by the endpoint $q=\infty$, so the bound of
Theorem~\ref{thm:Lq-rigidity}(iii) is also attained with $E_\infty=P$.

\begin{theorem}[Points on a geodesic]\label{thm:Lq-collinear}
Let $n\ge3$ and $k=\lfloor(n-1)/2\rfloor$.
\begin{enumerate}[label=\textup{(\roman*)}]
\item If $\mathcal A$ consists of $n$ points on a geodesic $\gamma$ and
the center path is nonconstant, then $E_q\in\gamma$ for every
$q\in(1,\infty]$, and $\nu(P)+\mathbf 1_{\{E_\infty=P\}}\le k$ for every
$P\in M$.
\item Conversely, for every geodesic $\gamma$, every $P\in\gamma$, every
$R>0$ and all $1<q_1<\cdots<q_k<\infty$, there are $n$ points on $\gamma$,
within distance $R$ of $P$, with $E_q=P$ exactly for
$q\in\{q_1,\ldots,q_k\}$. For $q_i=i+1$ and $R=1$ one can take one point at
distance $1$ from $P$ and, for $j=1,\ldots,k$, a point of multiplicity $2$
at distance $\sin^2\bigl(j\pi/(2k+2)\bigr)$, on alternating sides of $P$
with the outermost double point opposite to the single point, together
with the point $P$ when $n$ is even.
\end{enumerate}
\end{theorem}

\begin{proof}
(i) Geodesics are closed and convex, so $E_q\in\gamma$ by
Lemma~\ref{lem:Lq-basic}(iii). The rest is the proof of
Theorem~\ref{thm:collinear-multiplicity}, with \eqref{eq:tail-identity} and
Lemma~\ref{lem:laplace-zeros} in place of Proposition~\ref{prop:hyp-tail}
and Lemma~\ref{lem:Lq-basic}(ii) in place of
Proposition~\ref{prop:equilibrium}.

(ii) The construction is that in the proof of
Theorem~\ref{thm:prescribed}(i), for $\mathcal P$ consisting of distinct
finite parameters, with $k_i=r^{q_i-2}$ and \eqref{eq:tail-identity} in place of
\eqref{eq:hyp-tail}, and with the points $\exp_P(\pm t\boldsymbol u)$ on
$\gamma$, where $\boldsymbol u$ is a unit tangent vector of $\gamma$ at
$P$. For $q_i=i+1$, $\varepsilon$ can be chosen explicitly.
Substitute $r=\sin^2(\theta/2)$, so that
$r^{i-1}\,dr=2^{-i}(1-\cos\theta)^{i-1}\sin\theta\,d\theta$ for
$1\le i\le k$. Each product $(1-\cos\theta)^{i-1}\sin\theta$ is a linear
combination of $\sin(l\theta)$ with $1\le l\le i$. The function
$\pm\operatorname{sgn}\sin((k+1)\theta)$ has a Fourier sine series on
$(0,\pi)$ containing only the frequencies $(k+1)(2m+1)$, so it is
orthogonal to all of them, and it changes sign exactly at
$\theta=j\pi/(k+1)$, that is, at $r=\sin^2(j\pi/(2k+2))$. It satisfies the
conclusion of Lemma~\ref{lem:canonical}, so it can serve as $\varepsilon$.
\end{proof}

\noindent The distances $\sin^2(j\pi/(2k+2))$ are the zeros of the
Chebyshev polynomial of the second kind $U_k(1-2r)$, the canonical points
of Markov's theorem on best $L^1$ approximation \cite{PinkusBook}.

\begin{example*}
Identify a geodesic with $\R$, take $P=0$, and let
\[
 \mathcal A=\{-\tfrac34,-\tfrac34,\tfrac14,\tfrac14,1\},
\]
the case $n=5$, $k=2$ of Theorem~\ref{thm:Lq-collinear}(ii). The scalar force at $0$ is
\[
 f(q)=1+2\bigl(\tfrac14\bigr)^{q-1}-2\bigl(\tfrac34\bigr)^{q-1},
\]
so $f(2)=f(3)=0$, while $E_\infty=\frac18$. The family is therefore
nonconstant, and Theorem~\ref{thm:Lq-collinear}(i) shows that $E_q=0$
exactly for $q\in\{2,3\}$: the mean and the $L^3$ center coincide,
although the configuration is not symmetric about $0$.
\end{example*}

\section{Further questions}\label{sec:questions}

We have bounded how often the center map can take a prescribed value
and constructed configurations that attain the bounds. We now ask about
its stationary parameters, the radial weights for which the tail bound
holds, and the cases in which the cut radii can be found explicitly.

\emph{Stationary parameters.}
For a configuration of $n$ points in $\Hh^2$, at how many finite
parameters can a nonconstant center map satisfy $\dot C_p=0$?
Let $B(n)$ be the supremum of this number over all such configurations.

Corollary~\ref{cor:three-points} gives $B(3)=0$. For $n\ge4$, we have
$B(n)\ge n-3$. Indeed, fix a point $P$ and choose $n-2$ distinct finite parameters.
The proof of Theorem~\ref{thm:prescribed}(ii) constructs a
configuration for which $C_p=P$ precisely at those parameters. The inner points lie on a geodesic through $P$, and
reflection in this geodesic exchanges the two outer points. By
Proposition~\ref{prop:equivariance} and uniqueness, every center lies
on the geodesic. Its signed distance from $P$ vanishes at each
prescribed parameter, so Rolle's theorem yields a stationary parameter
strictly between each consecutive pair.

For each $n$, $B(n)$ is finite. In upper-half-plane coordinates, the
equilibrium equation and the stationarity condition
\eqref{eq:velocity-center} define the stationary parameters in the
real exponential field, with the configuration as parameter. For a
nonconstant analytic path these parameters contain no interval.
O-minimality therefore makes each such set finite \cite{Wilkie}, and
uniform finiteness gives a bound depending only on $n$
\cite[Chapter~3]{vdDries}. What is the sharp bound? Already for four
points, can a nonconstant center map be stationary at two distinct
parameters?

\emph{Other radial weights.}
The shell argument applies to weights of the form
$w_p(r)=\tau(r)e^{p\lambda(r)}$; the tail bound requires a further
property of their radial derivatives. Suppose that $x=\lambda(r)$ is
a smooth increasing change of coordinate, and let
$a(x)=\tau(\lambda^{-1}(x))$. In this coordinate the tail kernel is
\[
 dw_p=e^{px}\bigl(a'(x)+p\,a(x)\bigr)\,dx.
\]
For the hyperbolic $p$-barycenters, $a(x)=\sqrt{1-e^{-2x}}$ and
$p\ge0$; for the $L^q$ centers, $a\equiv1$ and $p=q-1>0$.
In both cases, if $h$ is a finite step function that vanishes for all
sufficiently large $x$ and is nonzero on an interval, its transform has
at most $S^-(h)$ zeros, counted with multiplicity
(Theorem~\ref{thm:kernel-vd} and Lemma~\ref{lem:laplace-zeros}).
For which other positive functions $a$, and on which parameter
intervals, does this bound hold? The proof of
Lemma~\ref{lem:dirichlet} suggests looking for sufficient conditions
on $1/a$.

\emph{Explicit cut radii.}
For a fixed outer radius $R$ and a prescribed multiset of finite
parameters, Lemma~\ref{lem:canonical} determines the cut radii in
the construction of Theorem~\ref{thm:prescribed} uniquely. For the
parameters $1,3,\ldots,2m-1$, Remark~\ref{rem:odd-explicit} gives
these radii in terms of the zeros of a Chebyshev polynomial. Is there
a comparable formula, or a recurrence for the radii, for the list
$0,1,\ldots,m-1$? Its first two parameters give Julia's zero and the
hyperbolic zero.

\section*{Statements and declarations}

\noindent\textbf{Funding.} No funding was received.

\smallskip
\noindent\textbf{Competing interests.} The author has no competing interests.

\smallskip
\noindent\textbf{Data availability.} Not applicable.

\end{document}